\documentclass[a4paper,12pt,reqno]{amsart}
\usepackage{amssymb}
\usepackage[dvips]{graphicx}
\usepackage[all]{xy}
\usepackage{url}
\usepackage{hyperref}

\input epsf.tex
\newdimen\xsize
\newdimen\oldbaselineskip
\newdimen\oldlineskiplimit
\def\restorelineskip{\baselineskip=\oldbaselineskip%
\lineskiplimit=\oldlineskiplimit}
\def\putm[#1][#2]#3{
\hbox{\vbox to 0pt{\parindent=0pt%
\vskip#2\xsize\hbox to0pt{\hskip#1\xsize $#3$\hss}\vss}}}%
\long\def\Line#1{\hbox to \hsize{#1}}
\def\putt[#1][#2]#3{
\vbox to 0pt{\noindent\hskip#1\xsize\lower#2\xsize%
\vtop{\restorelineskip#3}\vss}}

\makeatletter
\def\xbig[#1]#2{{\hbox{$\m@th\left#2\vbox to#1\xsize{}%
\right.\n@space$}}}
\makeatother
\def\xlar[#1]#2{%
\smash{\mathop{ \hbox to #1\xsize{\leftarrowfill}}\limits^{#2}}}
\def\xrar[#1]#2{%
\smash{\mathop{ \hbox to #1\xsize{\rightarrowfill}}\limits^{#2}}}
\def\xline[#1]{\hbox to #1\xsize{\leaders\hrule\hfill}}

\DeclareFontFamily{U}{rsf}{\skewchar\font'177}%
\DeclareFontShape{U}{rsf}{m}{n}{<-6>rsfs5<6-8>rsfs7<8->rsfs10}{}%
\DeclareFontShape{U}{rsf}{b}{n}{<-6>rsfs5<6-8>rsfs7<8->rsfs10}{}%
\DeclareMathAlphabet\RSFS{U}{rsf}{m}{n}
\SetMathAlphabet\RSFS{bold}{U}{rsf}{b}{n}
\def\sf#1{{\mathsf{#1}}}

\def\slsf{\slshape \sffamily }

\def\msmall#1{\mathchoice{\hbox{\small$\displaystyle {#1}$}}{#1}{#1}{#1}}

\let\xrar=\xrightarrow

\def\bb{{\mathbb B}}

\def\cc{{\mathbb C}}
\def\dd{{\mathbb D}}

\def\rr{{\mathbb R}}

\def\nn{{\mathbb N}}

\def\pp{{\mathbb P}}

\def\arg{{\sf{Arg}}}

\def\dim{\sf{dim}\,}

\def\dist{\sf{dist}\,}

\def\ext{\sf{ext}}

\def\id{\sf{Id}}

\def\i{\sf{i}}

\def\lim{\mathop{\sf{lim}}}

\def\max{\sf{max}}
\def\min{\sf{min}}

\def\ret{\sf{r}}

\def\Reg{\sf{Reg}\,}

\def\Sing{\sf{Sing}\,}

\def\supp{\sf{supp}\,}
\def\sup{\sf{sup}\,}

\def\n{{\mathrm{n}}}

\def\to{\rightarrow}

\def\eps{\varepsilon}

\def\<{\langle}\let\la=\<
\def\>{\rangle}\let\ra=\>
 \let\bs=\bss 
\def\comp{\Subset}
\def\d{\partial}

\def\ddef{\mathrel{{=}\raise0.3pt\hbox{:}}}
\def\deff{\mathrel{\raise0.3pt\hbox{\rm:}{=}}}

\def\fraction#1/#2{\mathchoice{{\msmall{ #1\over#2}}}%
{{ #1\over #2 }}{{#1/#2}}{{#1/#2}}}
\def\norm#1{\left\Vert{#1}\right\Vert}

\def\le{\leqslant}

\def\emptyset{\varnothing}
\def\scirc{\mathop{\mathchoice{\hbox{\small$\circ$}}{\hbox{\small$\circ$}}%
{{\scriptscriptstyle\circ}}{{\scriptscriptstyle\circ}}}}
\def\longpoints{\leaders\hbox to 0.5em{\hss.\hss}\hfill \hskip0pt}
\def\stateskip{\smallskip}
\def\state#1. {\stateskip\noindent{\bf#1. }} 
\def\statep#1. {\stateskip\noindent{\bf#1 }} 
\def\proof{\state Proof. }

\def\Chi{\raise 2pt\hbox{$\chi$}}

\def\ie{\hskip1pt plus1pt{\sl i.e.\/,\ \hskip1pt plus1pt}}

\def\sli{{\sl i)} } 
\def\slii{{\sl i$\!$i)} } 
\def\sliii{{\sl i$\!$i$\!$i)} }
\def\sliv{{\sl i$\!$v)} }

\def\star{\mathop{\msmall{*}}}
\def\barr#1{\mskip1mu\overline{\mskip-1mu{#1}\mskip-1mu}\mskip1mu}
\def\Chi{\raise 2pt\hbox{$\chi$}}

\let\phI=\phi\let\phi=\varphi\let\varphi=\phI
\let\cal=\mathcal

\def\calc{{\cal C}}

\def\calk{{\cal K}}

\def\calo{{\cal O}}

\def\trans{\pitchfork}
\def\eps{\varepsilon}

\def\bs{\backslash}

\def\comp{\Subset}
\def\d{\partial}

\def\1{{1\mkern-5mu{\rom l}}}

\def\ge{\geqslant}

\def\fraction#1/#2{\mathchoice{{\msmall{ #1\over#2}}}%
{{ #1\over #2 }}{{#1/#2}}{{#1/#2}}}

\def\le{\leqslant}

\def\emptyset{\varnothing}

\def\qed{\ \ \hfill\hbox to .1pt{}\hfill\hbox to .1pt{}\hfill $\square$\par}

\def\comment#1\endcomment{}

 \belowdisplayskip=\abovedisplayskip

\def\lineeqqno(#1){\hfill\llap{\vbox to 10pt%
{\vss\begin{align} \eqqno(#1)\end{align}\vss}}\vskip1pt}

\advance\headheight 1.2pt

\def\ShowwLLabel#1{}

\def\thechpt{\Roman{chpt}}

\def\newchapt[#1]#2{\newpage%
\refstepcounter{chpt}\setcounter{subsection}{0}%
\setcounter{thm}{0}\setcounter{defi}{0}%
\setcounter{rema}{0}\setcounter{exrc}{0}%
\renewcommand{\thesubsection}{\thechpt.\arabic{subsection}}%
\section*{\begin{center}\huge \bf Chapter \thechpt\\
#2 \end{center}}\label{#1}%
\ \smallskip%
\markboth{Chapter \thechpt}{#2}%
}
\def\newsect[#1]#2{\refstepcounter{section}\setcounter{equation}{0}%
\renewcommand{\thesubsection}{\arabic{section}.\arabic{subsection}}%
\section*{\arabic{section}.
#2}\vspace{-20pt}\label{#1}\vspace{20pt}%
\markboth{Section \arabic{section}}{#2}}

\def\newlect[#1]#2{\refstepcounter{section}%
\renewcommand{\thesubsection}{\arabic{section}.\arabic{subsection}}%
\section*{Lecture \arabic{section}\\
#2}\label{#1}%
\markboth{Lecture \arabic{section}}{#2}}

\def\newprg[#1]#2{\refstepcounter{subsection}%
\subsection*{{\thesubsection.\ #2}} \label{#1}%
}

\def\newappx[#1]#2{%
\refstepcounter{appx}\setcounter{section}{0}%
\renewcommand{\thesubsection}{A\arabic{appx}.\arabic{subsection}}%
\section*{Appendix \arabic{appx}\\ #2}
\label{#1}%
\markboth{Appendix A\arabic{appx}}{#2}
}

\newtheorem{asser}{Assertion}
   \def\newasser#1{\begin{asser} \label{#1}}

\newtheorem{thm}{Theorem}[section]
   \def\newthm#1{\begin{thm}\label{#1}}

\newtheorem{nnthm}{Theorem}
   \def\newthm#1{\begin{nnthm}\label{#1}}

\newtheorem{lem}{Lemma}[section]
   \def\newlemma#1{\begin{lem} \label{#1}}

\newtheorem{prop}{Proposition}[section]
   \def\newprop#1{\begin{prop}\label{#1}}

\newtheorem{nnprop}{Proposition}
   \def\newprop#1{\begin{nnprop}\label{#1}}

\newtheorem{claim}{Claim}[section]
   \def\newclaim#1{\begin{claim} \label{#1}}

\newtheorem{corol}{Corollary}[section]
   \def\newcorol#1{\begin{corol} \label{#1}}

\newtheorem{nncorol}{Corollary}
   \def\newcorol#1{\begin{nncorol} \label{#1}}

\newtheorem{defi}{Definition}[section]
   \def\newdefi#1{\begin{defi} \label{#1}\rm }

\newtheorem{nndefi}{Definition}
   \def\newdefi#1{\begin{nndefi} \label{#1}\rm }

\newtheorem{exmp}{Example}[section]
   \def\newexmp#1{\begin{exmp} \label{#1}\rm }

\newtheorem{nnexmp}{Example}
   \def\newexmp#1{\begin{nnexmp} \label{#1}\rm }

\newtheorem{problem}{Problem}
   \def\newexmp#1{\begin{problem} \label{#1}\rm }

\newtheorem{exrc}{Exercise}
   \def\newexrc#1{\begin{exrc} \label{#1}\rm }

\newtheorem{rema}{Remark}[section]
   \def\newrema#1{\begin{rema} \label{#1}\rm }

\newtheorem{nnrema}{Remark}
   \def\newrema#1{\begin{nnrema} \label{#1}\rm }

\def\eqqno(#1){\label{(#1)}}
\def\eqqref(#1){(\ref{(#1)})}

\title{Discrete and Continuous Versions of \\ the Continuity Principle}
\author{S. Ivashkovich}
\date{\today}

\address{
Universit\'e de Lille-1, UFR de Math\'ematiques, 59655 Villeneuve
d'Ascq, France} \email{serge.ivashkovych@univ-lille.fr}

\subjclass[2000]{Primary - 32D15, Secondary - 32D10}
\keywords{Continuity Principle, holomorphic function, meromorphic function, analytic set, Hausdorff metric.}

\begin{document}
\begin{abstract}
The goal of this paper is to  present a certain generalization of the classical  Kontinuit\"atssatz 
of Behnke for holomorphic/meromorphic functions in terms of the lift to the envelope of holomorphy. 
We consider two non-equivalent formulations: ``discrete'' and ``continuous'' ones. Giving a proof 
of the ``discrete'' version we, somehow unexpectedly, construct a counterexample to the ``continuous'' 
one when convergence/continuity of analytic sets is considered in Hausdorff topology or, even in the 
stronger topology of currents. But we prove the ``continuous'' version of the Kontinuit\"atssatz if 
continuity is understood with respect to the Gromov topology. Our formulations seem to be not yet 
existing in the literature. A number of relevant examples and open questions is given as well.
\end{abstract}

\maketitle

\setcounter{tocdepth}{1}
\tableofcontents

\newsect[INT]{Introduction}
\newprg[INT.d]{Discrete Version}
Let $X$ be a complex manifold. By an analytic set with boundary in $X$ we mean an analytic set $C$ in some 
open subset $U\subset X$ and we define the boundary of $C$ as  $\d C\deff \overline{C}\cap \d U$. If $U$ is 
relatively compact in $X$ we say that $(C,\d C)$ is a {\slsf compact} analytic set with boundary in $X$.
We say that compact analytic sets with boundary $(C_k,\d C_k)$ converge to a compact analytic set with boundary 
$(C_0,\d C_0)$ in {\slsf Hausdorff} topology if both $\overline{C}_k\to \overline{C}_0$ and $\d C_k\to \d C_0$ 
in Hausdorff metric, see Definition \ref{hausd-top} in Section \ref{BEHNKE}. Notice that we do not require 
that $C_k$ and $C_0$ are analytic sets in some common open $U\subset X$. It may well happen that every 
$C_k$ is analytic in  its own $U_k\comp X$, as well as $C_0$ is analytic in some $U_0\comp X$, 
and all $U_k$ are distinct. All analytic sets in this paper 
are supposed, if the opposite is not explicitly stated, to be {\slsf proper}, \ie $\dim C \le \dim X-1$ and
having all their irreducible components of dimension $\ge 1$. By ``components'' we mean irreducible components 
of $C$ and not of $\d C$, even if the latter has some analytic structure. Our first goal in this paper is to 
prove the following ``discrete'' version of the Continuity Principle.

\newpage

\begin{nnthm} {\slsf (Continuity Principle - I).}
\label{cpd}
Let $D$ be a domain in a Stein manifold $X$ and let $\{(C_k,\d C_k)\}_{k\in \nn}$ be a sequence of pure $q$-dimensional 
compact analytic sets with boundary in $D$, $1\le q< n =\dim_{\cc}X$. Suppose that $(C_k,\d C_k)$  converge in Hausdorff 
topology to a pure $q$-dimensional compact analytic set with boundary $(C_0,\d C_0)$ in $X$ such that $\d C_0\comp D$. 
Then $C_0$ can be lifted to the envelope of holomorphy $(\widehat{D}, \pi)$ of $D$. Namely there exists a 
compact analytic set with boundary $(\widehat{C}_0, \d \widehat{C}_0)$ in $\widehat{D}$ such that:

\sli the restriction  $\pi|_{\widehat{C}_0}:\widehat{C}_0\to C_0$ is proper surjective and one-to-one near the 
boundaries.

\slii $(\i (C_k), \i (\d C_k))$ converge to $(\widehat{C}_0, \d \widehat{C}_0)$ in Hausdorff topology, here 
$\i :D\to \widehat{D}$ is the 

\quad canonical inclusion.
\end{nnthm}

\noindent In particular, we shall see that  $\d\widehat{C}_0=\i (\d C_0)$, \ie  the boundary of 
$\widehat{C}_0$ is precisely $\i (\d C_0)$ and nothing more. The statement of Theorem \ref{cpd} means that 
every holomorphic function in $D$ extends as a holomorphic function to a {\slsf fixed}, \ie independent of 
a function, neighborhood of $\widehat{C}_0$ producing thus a sort of a {\slsf multivalued analytic extension}
of holomorphic in $D$ functions to a fixed neighborhood of $C_0$. 

\begin{nnrema} \rm
\label{mer}
a) Via the result of \cite{KS} the statement of Theorem \ref{cpd} gives also a multivalued analytic extension 
for meromorphic functions in $D$ to a neighborhood of $C_0$. 

\noindent b) Analytic extension in Theorem \ref{cpd} is not singlevalued in general when $q<n/2$, 
see Part-I of Example \ref{moy-exmp1} in Section \ref{LIFT-C}. But if $q\ge n/2$ the extension {\slsf is} singlevalued. 
We prove the following

\begin{nnprop}
\label{cp1}
If under the assumptions of Theorem  \ref{cpd} one supposes, in addition, that $q\ge n/2$ then the lift $\widehat{C}_0$ of $C_0$  
is singlevalued, \ie $\pi|_{\widehat{C}_0}: \widehat{C}_0\to C_0$ is  an isomorphism.
\end{nnprop}

\noindent c) The lift $\widehat{C}_0$ of $C_0$ in Theorem \ref{cpd} is constructed as follows. Using the fact that the canonical 
inclusion $\i :D\to \widehat{D}$ is a biholomorphism onto its image we can find a 
neighborhood $W_0$ of $\d C_0$ and a neighborhood $\widehat{W}_0\subset \widehat{D}$ of $\d\widehat{C}_0 
\deff \i (\d C_0)$ such that $\i: W_0\to \widehat{W}_0$ is a biholomorphism. We define then $\widehat{C}_0$ 
as the union of all irreducible components of $\pi^{-1}(C_0)$ which intersect $\widehat{W}_0$ and we 
prove that this $\widehat{C}_0$ satisfies the conclusion of Theorem \ref{cpd}.
\end{nnrema}

\newprg[INT.c]{Continuous Version}
It is natural to consider the following ``continuous'' version of the Continuity Principle. Let $D$ be a domain in a Stein 
manifold $X$ and let $\{(C_t,\d C_t)\}_{t\in [0,1]}$ be a continuous in Hausdorff topology family of pure $q$-dimensional 
compact analytic sets with boundary in $X$, $1\le q<n =\dim_{\cc}X$, such that $\bar C_0\subset D$ and $\d C_t\subset D$ for 
all $t\in [0,1]$. The question is: can {\it $\{(C_t, \d C_t)\}_{t\in [0,1]}$ be lifted to the envelope of holomorphy 
$(\hat D,\pi)$ of $D$?} Namely, does there exist a continuous in Hausdorff topology family 
$\{(\widehat{C}_t, \d\widehat{C}_t)\}_{t\in [0,1]}$ of pure $q$-dimensional compact analytic sets with boundary in $\widehat{D}$ 
such that: 

\smallskip\sli  for all $t\in [0,1]$ the restriction $\pi|_{\widehat{C}_t}:\widehat{C}_t \to C_t$ is proper surjective and one to 
one near

\quad  the boundaries $\d\widehat{C}_t=\i (\d C_t)$  of $\widehat{C}_t$ and $\d C_t$ of $C_t$;

\smallskip\slii $(\widehat{C}_t, \d \widehat{C}_t) = (\i (C_t), \i (\d C_t))$ for $t\in [0,1]$ close to zero.

As a candidate for the lifts of $C_t$ should be $\widehat{C}_t$ constructed in the same way as $\widehat{C}_0$ 
in Theorem \ref{cpd}, \ie $\widehat{C}_t$ should be the union of all irreducible components of $\pi^{-1}(C_t)$ 
which intersect $\widehat{W}_0$. This time $\widehat{W}_0$ is biholomorphically mapped by $\pi$ to an appropriate
neighborhood $W_0$ of $\bigcup_{t\in [0,1]}\d  C_t$. If dimension $q$ of $C_k$-s satisfies the bound $q\ge n/2$ as 
in Proposition \ref{cpd} this ``continuous'' version of the CP  holds true. It is essentially a corollary of the 
proofs of Theorem \ref{cpd} and Proposition \ref{cp1}.

\begin{nnprop} 
\label{cp2}
Let $D$ be a domain in a Stein manifold $X$ and let $\{(C_t,\d C_t)\}_{t\in [0,1]}$ be a continuous in Hausdorff 
topology family of compact analytic sets with boundary in $X$ of pure dimension $\dim X/2\le q \le \dim X-1$ 
such that $\bar C_0\subset D$ and $\d C_t\subset D$ for all $t\in [0,1]$. Then the family $(C_t,\d C_t)$ can be
continuously lifted to the envelope of holomorphy $(\hat D, \pi)$ of $D$. Moreover, if $(\widehat{C}_t, \d
\widehat{C}_t)$ is a corresponding lift then $\pi|_{\widehat{C}_t}:\widehat{C}_t\to C_t$ is an isomorphism 
for every $t\in [0,1]$.
\end{nnprop}
Let us state this proposition equivalently but somewhat differently.
\begin{nncorol}
\label{cp2-cor}
In the conditions of Proposition \ref{cp2} every holomorphic/meromorphic function $f\in D$ can be analytically 
continued along $\{(C_t, \d C_t)\}_{t\in [0,1]}$. Namely, there exists a family $\{f_t\}_{t\in [0,1]}$ of 
holomorphic/meromorphic functions in $V_t$, where $V_t$ is a neighborhood of $\bar C_t$, such that:

\sli $f_0=f|_{V_0}$, where $V_0$ is a neighborhood of $C_0$ contained in $D$;

\slii $f_{t_1}=f_{t_2}$ on $V_{t_1}\cap V_{t_2}$ for $t_1$ close to $t_2$.
\end{nncorol}

We do not claim, and this is not true in general, that $f_t$-s glue together to a singlevalued function in 
some subdomain of $X$ bigger than $D$. 
\begin{nnrema} \rm
It is worth of noticing that when $\dim X=2$ (and therefore $C_t$ are curves) the ``continuous'' version of CP holds true
in the form of Proposition \ref{cp2} and Corollary \ref{cp2-cor}.
\end{nnrema}

Somewhat surprisingly the statement of the the ``continuous'' version of the Continuity Principle does not hold true if 
$q<\dim X/2$. In Section \ref{LIFT-C} we construct the following 

\begin{nnexmp}
\label{moy-exmp}
There exists a domain $D\subset \cc^3$ and a continuous family $\{(C_t,\d C_t)\}_{[-1, 1]}$ of complex curves 
with boundary such that all $C_t$ except $C_0$ are smoothly imbedded, $C_0$ is immersed with one double point 
and such that this family possesses the following properties:

\sli $\bigcup_{[-1, 1]} \d C_t\comp D$ and $\bar C_{-1}\comp D$;

\slii ${C_t}$ can be continuously lifted to $\widehat{D}$ up to $0$. 

\sliii For every $t>0$ the lift $\widehat{C}_t$ of $C_t$ is irreducible but has boundary components other 

\quad than  $\i (\d C_t)$. Also  
\[
\lim_{H}\widehat{C}_t\not= \widehat{C}_0 \qquad\text{ as } \qquad  t\searrow 0.
\]
\end{nnexmp}

Here $\lim\limits_H$ denotes the Hausdorff limit. The family $C_t$ of this example is continuous not only
in Hausdorff topology but also in a stronger topology of currents. But it is discontinuous at zero 
in Gromov topology, see more about all this in Section \ref{LIFT-C}. It turns out that the continuity in 
the latter topology is sufficient for the validity of the  ``continuous'' version of the CP. Recall that a 
compact complex curve with boundary over a complex manifold $X$ is a pair $(C,u)$, where $C$ is a compact
analytic space of dimension one with only nodes as singularities (smooth near the boundary) and $u:C\to X$ 
is a holomorphic mapping.





\begin{nnthm} {\slsf (Continuity Principle - II).}
\label{cpc}
Let $D$ be a domain in a Stein manifold $X$ and let $\{(C_t,u_t)\}_{t\in [0,1]}$ be a family of compact 
complex curves with boundary over $X$ which is continuous in Gromov topology. Suppose that $u_0(C_0)\subset D$ 
and $u_t(\d C_t)\subset D$ for all $t\in [0,1]$. Then this family can be lifted to the envelope of holomorphy 
$(\hat D,\pi)$ of $D$. Namely, there exist holomorphic mappings $\hat u_t:C_t\to \hat D$ such that:

\sli $\{(C_t,\hat u_t)\}_{t\in [0,1]}$ is continuous in Gromov topology family over $\hat D$;

\slii $\pi\circ \hat u_t = u_t$ for all $t\in [0,1]$.
\end{nnthm}
Let us formulate a simple particular case of this theorem which could be useful in applications.

\begin{nncorol}
\label{cp-c2}
Let $D$ be a domain in a Stein manifold $X$ and let $\{(C_t,u_t)\}_{t\in [0,1]}$ be a family of compact 
complex curves with boundary over $X$ which is continuous in $\calc^2$-topology. Suppose that $u_0(C_0)\subset D$ 
and $u_t(\d C_t)\subset D$ for all $t\in [0,1]$. Then this family can be lifted to the envelope of holomorphy 
$(\hat D,\pi)$ of $D$.
\end{nncorol}
If one supposes that $u_t:C_t\to X$ are imbeddings then this statement follows from the classical 
Behnke's Kontinuit\"atssatz, see discussion in section \ref{BEHNKE}. But in general in this Corollary 
the image $u_t(C_t)$ can have nodes and cusps, see Example \ref{r2-exmp} in Section \ref{LIFT-C}.

For a non-Stein $X$ the statements of  Theorems \ref{cpd} and \ref{cpc} doesn't hold true. At the end of 
Section \ref{LIFT-C} we give the following

\begin{nnexmp}
\label{nash-exmp}
There exist a sequence $(C_k, \d C_k) = (\phi_k(\Delta), \phi_k(\d \Delta))$ of imbedded analytic disks over a 
certain complex projective threefold $X$ converging in Gromov topology to a compact complex curve with boundary
$(C,\d C)$ such that:

\sli the limit $C$ is the union of an imbedded disk $\phi (\bar\Delta)$ and a rational curve $C_0$ (a bubble);

\slii there exists a holomorphic function in a domain $D\supset \bigcup_kC_k\cup \phi (\d\Delta)$ which has an 

\quad essential singularity (\ie is not even meromorphic) along $C_0$.

\sliii This sequence can be included to a continuous in Gromov topology family of complex 

\qquad curves with boundary $\{(C_t,\d C_t)\}_{t\in [0,1]}$ (by setting $C_k = C_{\frac{1}{k}}$) such that 

\quad $\bigcup_{t\in [0,1]}\d C_t\comp D$ and $\bigcup_{t\in (0,1]} C_t\subset D$.
\end{nnexmp}

\begin{nnrema} \rm 
To our best knowledge the statements of Theorems \ref{cpd} and \ref{cpc} are more general that  the existing
versions of the Continuity Principle (in what follows CP for short) in the literature. The closest one we know 
about is the result of Chirka and Stout in \cite{CS}, where more is assumed about convergence of $C_k$, they 
should converge in the topology of currents. We say more about the result of \cite{CS} in the discussion after
Example \ref{moy-exmp1} in Section \ref{LIFT-C}, see Remark \ref{ch-st}, where we point out a problem in the 
approach of \cite{CS}. And this problem is connected to the failure of the ``continuous'' version of CP.
\end{nnrema}

\smallskip\noindent{\slsf The structure of the paper.} {\slsf 1.} In section \ref{BEHNKE} we shall prove the CP 
in the form of Behnke replacing $\calc^2$-convergence by the Hausdorff one. We also give a version of the 
Behnke-Sommer result. The distinguished feature of these formulations is that one supposes that the limit $C_0$ 
is a smooth manifold. In that case the extension is single-valued. We also formulate some open questions.

\smallskip\noindent{\slsf 2.} In section \ref{LIFT} we prove Theorem \ref{cpd}. The key point in the
proof of these statements is the ``lift of paths'' Lemma \ref{path-lift1} of section \ref{LIFT}. We prove 
there also Propositions \ref{cp1} and \ref{cp2}.

\smallskip\noindent{\slsf 3.} In section \ref{LIFT-C} we construct the Example \ref{moy-exmp}, discuss the 
approach of \cite{CS} and formulate one more open question, see Remark \ref{ch-st}. After that we recall 
the notions connected with the Gromov topology and prove Theorem \ref{cpc}. We also construct Example 
\ref{nash-exmp} there. At the end we formulate one more problem concerning a CP over non-Stein complex surfaces.

\smallskip\noindent{\slsf 4.} Along this paper we freely use the notions and results connected with envelopes 
of holomorphy. As a sources we recommend the first chapter of \cite{GR} and more recent exposition in \cite{Jr}.

\newsect[BEHNKE]{Continuity Principle in the form of Behnke}

To fix the notations and for the sake of future references we give here a version of the Behnke's Continuity 
Principle in the form a bit more general than it can be usually found in the literature but still different 
from Theorem \ref{cpd}. Recall the following notion. 

\begin{defi}
\label{hausd-top}
For compacts $A,B$ in a metric space $(X,d)$ the {\slsf Hausdorff distance} between $A$ and $B$ is defined as
\begin{equation}
\dist_H(A,B) \deff \inf\{\eps >0: B^{\eps} \supset A, A^{\eps}\supset B\}.
\end{equation}
\end{defi}

\noindent Here for  $\eps >0$ the set $A^{\eps}$ is called the {\slsf $\eps$-neighborhood} of $A$ and is defined as
\begin{equation}
\eqqno(eps-neib)
A^{\eps} \deff \{x\in X: d(x,A)<\eps\}.
\end{equation}
According to this definition a sequence $\{A_k\}$ of compact subsets of a metric space $(X,d)$ converges to 
a compact subset $A\comp X$ in {\slsf Hausdorff topology} if $\dist_H(A_k,A)$ tends to zero.

\begin{rema} \rm 
A more common notion of a distance between two compacts will be used in this paper as well:
\begin{equation}
\eqqno(dist-e)
\dist (A,B) \deff \inf\{ d(a,b):a\in A, b\in B\}.
\end{equation}
Convergence with respect to this distance will be not considered.
\end{rema}

Let $S$ be a finite disjoint union of smoothly imbedded circles in a complex manifold $X$.
By a {\slsf compact complex curve with boundary} $S$ in $X$ we understand a complex analytic subset 
$C$ of $X\setminus S$ of pure complex dimension one such that the union $C\cup S$ is compact in $X$ 
and near $S$ the set $\bar C$ is a smooth manifold with boundary $\d C =S$. By a {\slsf smooth} 
compact complex curve with boundary we understand a smooth compact complex submanifold with boundary 
of $X$ of complex dimension one.  In accordance with general 
notion from the Introduction we say that a sequence $(C_k,\d C_k)$ of compact complex curves with 
boundary converge to a compact complex curve with boundary $(C_0,\d C_0)$ in Hausdorff topology if both 
\[
\dist_H(\bar C_k, \bar C_0)\to 0 \quad\text{ and }\quad \dist_H(\d C_k, \d C_0)\to 0.
\]

\index{Theorem!Behnke} 

\begin{figure}[h]
\centering
\includegraphics[width=2.5in]{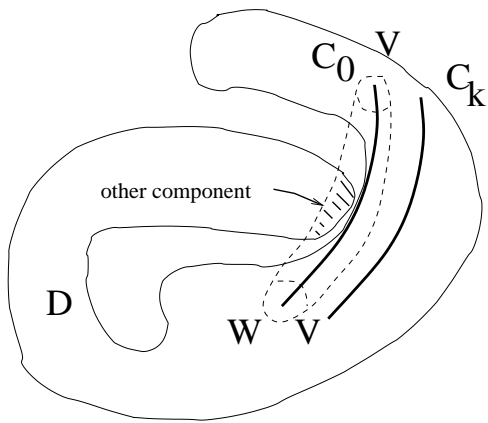}
\caption{Domain $D$ on this picture touches itself at the place indicated by an arrow. As a result an 
intersection of a neighborhood $W$ of $C_0$ with $D$ has two components and to the dashed component
$f$ might extend with different value.}
\label{other-comp}
\end{figure}

\begin{thm} 
\label{behnke}
Let $D$ be a domain in a complex manifold $X$, $\dim X= n\ge 2$, and let $\{(C_k, \d C_k)\}_{k\in \nn}$
be a sequence of compact complex curves in $D$ with boundary converging in $X$ to a {\slsf smooth compact} 
complex curve  with boundary $(C_0, \d C_0)$ in Hausdorff topology. Suppose that:

\sli $\d C_0\comp D$,

\slii $C_0$ has no irreducible components without boundary.

\smallskip\noindent Then there exist tubular neighborhoods $D\supset V \supset \d C_0$ and $W\supset \bar C_0$ such 
that for every holomorphic/meromorphic function $f$ in $D$ there exists a holomorphic/meromor\-phic function 
$\tilde f$ in $W$ with $\tilde f|_V = f|_V$. 
\end{thm}
\proof In other words for every $f\in \calo (D)$ its restriction $f|_V$ extends to a neighborhood $W$ of $\bar C_0$ 
and this $W$ doesn't depend on $f$. One should notice that $W\cap D$ might have connected components other than that 
which contains $V$, see Picture \ref{other-comp},  and $\tilde f|_{W\cap D}$ may not coincide with $f|_{W\cap D}$ 
on these components in general. 

By tubular neighborhoods in this theorem we mean an $\eps$-neighborhoods with 
respect to some properly chosen metric: 
\begin{equation}
\eqqno(tub-neib1)
V =  \d C_0^{\eps} = \{ x\in X: \dist (x, \d C_0) < \eps \} \quad \text{and} \quad 
W = \bar C_0^{\eps }= \{ x\in X: \dist (x, \bar C_0) < \eps \}.
\end{equation}
\begin{rema}\rm
\label{metric1}
Let us make a precision about a possible choice of such metrics. If $X=\cc^n$ we take the metric $d_{\Delta^n}$
induced by the polydisk norm. When working with domains $(\hat D,\pi)$ over $\cc^n$ we equip $\hat D$ with the pull-back 
of the polydisk metric and therefore $\pi : \hat D\to \cc^n$ will be a local isometry. The distance $\dist (\cdot,\cdot)$
as in \eqqref(dist-e) induced by the polydisk norm and we shall call  also the {\slsf polydisk distance}. 
If $X$ is Stein imbedd $X$ to some $\cc^N$ and take the distance induces by the polydisk norm.
In the case of general $X$ remark that one can find a Stein neighborhood of $\bar C_0$, see \cite{Siu}, and again 
one can use the induced polydisk distance.
\end{rema}

\smallskip\noindent 
{\slsf Case 1. Suppose first that $X=\cc^n$ and that $f$ is holomorphic.} Take a point $p\in D$ and write 
the Taylor expansion
\begin{equation}
\eqqno(taylor1)
f(z) = \sum_{|m|=0}^{\infty} a_m(f,p)(z-p)^m
\end{equation}
of $f$ with center at $p$. For every $m\in \nn^n$ coefficient $a_m(f, p)$ is a holomorphic function of $p$
in $D$. Let $d>0$ be a number smaller than the polydisk distance $d_0$ from $\d C_0$ to $\d D$. 
Using the assumption that $C_0$ is imbedded take $d$ sufficiently small in order for 
\[
\bar C_0^d\deff \bigcup_{p\in C_0}\Delta^n(p,d)
\]
to be a tubular neighborhood of $\bar C_0$. Take $d_1=\frac{d_0-d}{2}$ and set $\overline{\d C_0^{d+d_1}} = 
\{z: d(z, \d C_0)\le d+d_1\}$. Furthermore set $M\deff M(f,\overline{\d C_0^{d+d_1}})\deff\max\{ |f(z)|: z\in 
\overline{\d C_0^{d+d_1}}\}$. Then Cauchy inequalities gave us $|a_m(f,p)|\le M/d^{|m|}$ for all $p\in 
\overline{\d C_0^{d_1}}$. Therefore by the maximum modulus principle we have
\begin{equation}
\eqqno(cauchy1)
|a_m(f,p)|\le \frac{M}{d^{|m|}}
\end{equation}
for all $p\in C_k$ and all $k$ big enough. More precisely for all $k$ such that $\dist_H (\bar C_k, \bar C_0) \le d_1$. 
This implies that \eqqref(cauchy1) holds true also for all $p\in C_0\cap D$. Let $V=\d C_0^d$ be a $d$-neighborhood of 
$\d C_0$. We see that  $f|_V$ extends holomorphically to the tubular $d$-neighborhood of $C_0\cap D$. Let $\Omega$ be the 
maximal open subset of $C_0$ such that $f$ holomorphically extends to the tubular $d$-neighborhood of 
$\Omega$. This $\Omega$ is open and contains $C_0\cap D$. If $\Omega$ is not the whole of $C_0$ then 
take some $p_0\in \d \Omega$ and some $p\in \Omega$ close to this $p_0$. By the argument above $f$ holomorphically
extends to $\Delta^n(p,d)$ and henceforth (since $p$ can be taken arbitrarily close to $p_0$) to 
$\Delta^n(p_0,d)$. This proves that $\Omega$ must be the whole of $C_0$.

\smallskip\noindent{\slsf Case 2. $X$ is a Stein manifold.} Imbedd it to $\cc^N$. If $\i : X\to \cc^N$ 
is this embedding we denote $\i (X)$ still by $X$. Fix a sufficiently small neighborhood $V\supset X$ such that 
there exists a holomorphic retraction $\ret : V\to X$, see \cite{GR} Chap. VIII, Theorem C8. Set
$D^{\ret}\deff \ret^{-1}(D)$. Now we are in $\cc^N$ and our holomorphic ($\ret$-invariant!) $f\circ \ret $ extends 
to a neighborhood of $C_0$ staying to be $\ret$-invariant. This gives an extension of $f$ itself. 

\smallskip\noindent{\slsf Case 3. $X$ is Stein and $f$ is meromorphic.} Then $f$ can be represented 
as $f=h/g$ with $h,g\in \calo (D)$, see \cite{KS}. Now the conclusion of the theorem follows from the holomorphic case. 

\smallskip\noindent{\slsf Case 4. $X$ is an arbitrary complex manifold.} Remark that due to the assumptions of our 
theorem $C_0$ is Stein. Take a Stein neighborhood $X_1$ of $C_0$, see \cite{Siu} or \cite{De}, and reduce the situation 
to the Stein case. 

\smallskip\qed

\begin{rema} \rm
\label{be-som1}

\noindent{\bf a)} 
In \cite{Bh} this theorem was proved for $X=\cc^n$ assuming that complex manifolds with boundary $C_k$ (not 
necessarily curves) converge to $C_0$ in $\calc^2$-topology. In this case, \ie if $X=\cc^n$, complex manifolds/curves 
$C_k$ can be replaced by smooth real submanifolds of $\cc^n$ with boundary satisfying the maximum modulus 
principle in the sense that for every holomorhic function $f$ in a neighborhood of $\overline{C}_k$ one has
\begin{equation}
\label{max-pr}
\max_{z\in C_k} |f(z)| \le \max_{z\in \d C_k} |f(z)|.
\end{equation}
This generalization is due to Behnke-Sommer, see \cite{BS}. Let us state a version of it replacing $\calc^2$-convergence
by the Hausdorff one.

\begin{prop} 
\label{behn-pr1}
Let $D$ be a domain in a Stein manifold $X$, $\dim X= n\ge 2$, and let $\{(C_k, \d C_k)\}_{k\in \nn}$
be a sequence of compact real submanifolds with boundary in $D$, satisfying the maximum modulus principle,
and converging to an imbedded compact submanifold with boundary $(C_0, \d C_0)$ of $X$ in Hausdorff topology 
such that $\d C_0\comp D$. Then every holomorphic/meromorphic function $f$ in $D$ extends to a 
holomorphic/meromor\-phic function $\tilde f$ in a neighborhood of $\overline{C}_0$ which doesn't depend
on $f$.
\end{prop}
\end{rema}
As in the case of analytic sets Hausdorff convergence of $(C_k, \d C_k)$ to $(C_0,\d C_0)$ means that 
$\bar C_k\to \bar C_0$ and $\d C_k\to \d C_0$ in Hausdorff metric. The proof is literally the same as that 
of Theorem \ref{behnke}. 

\begin{problem} \rm Can one characterize compact submanifolds with 
boundary in $\cc^n$ satisfying the {\slsf maximum modulus principle}? Do they necessarily contain 
germs of complex curves through {\slsf any} point? It is not difficult to prove that this condition is 
sufficient.

\begin{prop}
\label{behn-pr2}
Let $C$ be a compact real submanifold of $\cc^n$ with boundary such that for every point $p\in C\setminus \d C$
there exist a germ of a non constant complex curve passing through $p$ which is contained in $C\setminus \d C$. 
Then $C$ satisfies the maximum modulus principle.
\end{prop}
\proof Suppose not, \ie there exists a holomorphic in a neighborhood of $C$ function $f$ and a point $p\in C
\setminus \d C$ such that $|f(p)| > \max_{z\in \d C}|f(z)|$. Set $K =\{ z\in C:f(z) =f(p)\}$. $K$ is a 
non-empty compact disjoint from $\d C$ and such that 

\medskip $\star$ $K$ contains a non-constant germ of a complex curve through each of its points. 

\medskip\noindent Take now the function $f_1(z) = z_1$ and let $p_1\in K$ be a point at which $|z_1|$ achieves 
its maximum. Set $K_1\deff\{ z\in K: f_1(z) = f_1(p_1)\}$. Then compact $K_1\comp K$ is non-empty and possesses 
the property $\star$. Pass to $f_2(z)= z_2$ and so on. Get a sequence $K\supset K_1\supset ...\supset K_n$ of 
non-empty compacts possessing property $\star$. Set $K_0\deff \cap_jK_j$ and get a contradiction: $K_0$ is 
non-empty, satisfies $\star$, all $z_1,...,z_n$ are constant on it, \ie $K_0$ is a singleton. This is impossible 
since it contains complex curves.

\smallskip\qed

The problem is: if the condition about complex curves is also necessary? Probably an easier question would be to 
characterize submanifolds in $\cc^n$ satisfying the {\slsf local} maximum modulus principle.
\end{problem}

\begin{problem} \rm 
The following open question related to the previous one is known for quite a long time as the {\slsf Problem of 
Rossi}. Let $M$ be a real analytic submanifold of $\cc^n$ such that the Levi form of $M$ degenerates at every 
point. Prove that for every point $p\in M$ there exists a germ of a non-constant complex curve through $p$ which 
is contained in $M$.
\end{problem}

\begin{rema} \rm
\label{be-so-re}
There is some discussion around these problems in \cite{BS} for pseudoconvex (but not strictly) 
$M$-s.
\end{rema}

\medskip Condition that $C_0$ is smooth in  Theorem \ref{behnke} is essential, \ie if $C_0$ is not imbedded then 
the extension to its neighborhood could be only multivalued as it is stated in Theorem \ref{cpd}. This will be 
illustrated by Example \ref{moy-exmp1} quoted in the Introduction.

\newsect[LIFT]{Lift to the envelope of holomorphy I: discrete case}

Now  let us give the proof of Theorem \ref{cpd} from the Introduction. It will be done in a number of steps. 
We start with the case $X=\cc^n$ first. Let us remark that compact analytic sets with boundary with all 
components of pure dimension $\ge 1$ do satisfy the maximum modulus principle. Indeed, let $C_0$ be a compact
analytic set with boundary in $\cc^n$, \ie $C_0$ is analytic in some relatively compact open $U_0$ and 
$\d C_0\deff \bar C\cap \d U_0$. Suppose that for some $p\in C_0$ and some $f$ holomorphic in a neighborhood 
of $\bar C_0$ one has $|f(p)| >\max\{|f(z)|: z\in \d C_0\}$. Then by continuity of $f$ one finds $U_1\comp 
U_0$ with smooth boundary sufficiently close to $\d U_0$ such that $|f(p)| >\max\{|f(z)|: z\in \bar C\cap 
\d U_1\}$ still holds. Contradiction with the usual maximum modulus principle.

\smallskip\noindent{\slsf Step 1. Local extension.} Let us state it in the form of a proposition.

\begin{prop}
\label{behn-pr3}
Let $D$ be a domain in $\cc^n, n\ge 2$, and let $\{(C_k,\d C_k)\}_{k\in \nn}$ be a sequence of proper compact 
analytic sets 
with boundary in $D$ with all irreducible components of positive dimension, converging in Hausdorff topology 
to a proper compact analytic set with boundary $(C_0,\d C_0)$ in $\cc^n$ such that $\d C_0\comp D$. Then for 
every point $p\in C_0$ and every $d>0$ smaller than the polydisk distance $d_0$ from $\d C_0$ to $\d D$
there exists a connected component $V$ of $\Delta^n(p,d)\cap D$ such that for any holomorphic in $D$ function 
$f$ the restriction $f|_V$ holomorphically extends to $\Delta^n(p,d)$. 
\end{prop}
\proof Take some $0<d_1<\min\{d,\frac{d_0-d}{4}\}$. Let $k$ be such that 
\[
\dist_H(\bar C_k, \bar C_0) + \dist_H (\d C_k , \d C_0) < d_1 ,
\]
and therefore 
\begin{equation}
\eqqno(dist-ck)
\dist (\d C_k, \d D) > d_0-d_1 > d + 3d_1.
\end{equation}
Fix a point $p\in C_0$  and take a point $p_{k}$ on $C_{k}$ on the polydisk distance $ < d_1$ from $p$. 
Let $O_k$ be the connected component of $\Delta^n(p,d_1)\cap D$ which contains $p_{k}$. Since compact analytic sets with 
boundary do satisfy the maximum modulus principle we get the estimate 
\begin{equation}
\eqqno(cauchy2)
|a_m(f,p_{k})|\le \frac{M(f, \overline{\d C_0^{d+2d_1}})}{(d+d_1)^{|m|}}
\end{equation}
for all $m\in \nn^n$ and every $f\in \calo (D)$ exactly in the same manner as in \eqqref(cauchy1). Here 
$\overline{\d C_0^{d+2d_1}}$ stands for  the closure of the $(d+2d_1)$-neighborhood of $\d C_0$. Remark that 
$\overline{\d C_0^{d+2d_1}}\comp D$ due to our choice of $d_1$. Therefore $f|_{O_k}$ holomorphically extends to 
$\Delta^n(p_{k},d + d_1)$. Since $\dist (p,p_{k})< d_1$ we have that $\Delta^n(p_{k},d+d_1)\supset \Delta^n(p,d)$ 
and therefore $f|_{O_k}$ extends holomorphically to $\Delta^n(p,d)$. Take a connected 
component $V$ of $\Delta^n(p,d)\cap D$ which contains $p_{k}$. Extension of $f|_{O_k}$ will be the extension of $f|_V$ 
as well. Proposition is proved.

\smallskip\qed

\begin{rema} \rm 
\label{proj-r}
{\bf a)} From this proposition we conclude that $\pi (\hat D)\supset \bar C_0$. Moreover, since $0<d<d_0$ can be taken 
arbitrarily close to $d_0$, we see that 
\begin{equation}
\eqqno(proj-f1)
\pi (\hat D)\supset \bar C_0^{d_0} \quad\text{ where } \quad d_0 = \dist (\d C_0 , \d D). 
\end{equation}
This follows from the estimate \eqqref(cauchy2) and Rossi's description of the envelope of holomorphy as the space 
of continuous homomorphisms from $\calo (D)$ to $\cc$, see \cite{Ro} or Chapter 1 of \cite{GR}. Indeed,
it is sufficient to prove that $\pi(\hat D)\supset \Delta^n(p_{k}, d+d_1)$, where $d, d_1$ and $p_k$ are 
taken as in the proof of Proposition \ref{behn-pr3}. Take any point $q\in 
\Delta^n(p_{k}, d+d_1)$. For $f\in \calo (D)$ denote by $\hat f$ its holomorphic extension to 
$\Delta^n(p_{k}, d+d_1)$ as above. Since $r\deff |q-p_{k}| < d+d_1$ we get from \eqqref(cauchy2)
\begin{equation}
\eqqno(proj-f2)
|\hat f(q)| \le \sum_{|m|=0}^{\infty}|a_m(f,p_{k})|r^{|m|}\le M(f, \overline{\d C_0^{d+2d_1}})
\sum_{|m|=0}^{\infty}\frac{r^{|m|}}{(d+d_1)^{|m|}} = 
\end{equation}
\[
= M(f, \overline{\d C_0^{d+2d_1}}) \left(1-\frac{r}{d+d_1}\right)^{-n} = 
M(f, \overline{\d C_0^{d+2d_1}})\left(\frac{d+d_1}{d+d_1-r}\right)^n.
\]
This means that the homomorphism $\phi_q$ defined as $\phi_q:f\to \hat f(q)$ is continuous and defines
a point in $\hat D$ over $q$.

\smallskip\noindent{\bf b)}
As it will be shown in Example \ref{moy-exmp1} this $V$ is not unique in general and for different $V$-s  the 
extensions
might be  different. Moreover, a component $V$ of $\Delta^n(p,d)\cap D$ in this Proposition might be such that 
$p\not\in \bar V$ in general! We know only that $V$ contains a point $p_{k}\in C_{k}$ close to $p$.
\end{rema}

\noindent{\slsf Step 2.  Lift to the envelope near the boundary.} 
Take a connected component $V$ of $\Delta^n(p, d)\cap D$ as in  Proposition \ref{behn-pr3}, \ie cutted
by $C_{k}$. We have a monomorphism $\calo (D)\to \calo (\Delta^n(p,d))$, namely holomorphic extensions of 
restrictions $f|_V$ to $\Delta^n(p,d)$. Denote this monomorphism as $\ext (\cdot |_V)$, \ie $f\to 
\ext (f|_V)$. This monomorphism is continuous by the estimate \eqqref(cauchy2) or, better by \eqqref(proj-f2). 
Now we can define a homomorphism 
\begin{equation}
\eqqno(hom-pv)
\phi_{p,V}:f\to \ext (f|_V)(p)
\end{equation}
from $\calo (D)$ to $\cc$, which is continuous as well. 
\begin{rema} \rm
Let us underline that monomorphism $\phi_{p,V}$ is well defined for all $p\in C_0$. But in general it depends also 
on the component $V$ of $\Delta^n(p,d)\cap D$ and eventually gives us the full lift of $\overline{C}_0$ to 
$\widehat{D}$. 
\end{rema}
But if $p\in \d C_0$ or, even more, $p\in D$ is just close to $\d C_0$ then there exists {\slsf only one}
component $V$ of $\Delta^n(p,d)\cap D$ namely $V = \Delta^n(p,d)$ itself and extension will be tautologically
to $\Delta^n(p,d)\subset D$. Therefore the following definition 
\begin{equation}
\eqqno(1-lift)
\d\widehat{C}_0\deff \{\phi_{p,V}: p \in \d C_0, V = \Delta^n(p,d)\cap D\}
\end{equation}
is correct and $\pi $ is one-to-one between an appropriately taken neighborhoods $\widehat{W_0}\supset 
\d \widehat{C}_0$ and $W_0\supset \d C_0$ correspondingly. Moreover it maps $\d\widehat{C}_0$ bijectively
to $\d C_0$. Its inverse is the restriction to $W_0$ of the canonical inclusion $\i :D\to \widehat{D}$.
In another words $\i (p) = \phi_{p,V}$ for $p\in W_0$ with $V=\Delta^n(p,d)$.

\smallskip
Set $\widetilde{C}_0\deff \pi^{-1}(C_0)$.  Since $\pi : \widehat{D}\to \cc^n$ is locally biholomorphic the set
$\widetilde{C}_0$ is closed in $\widehat{D}\setminus \pi^{-1}(\d C_0)$ and moreover $\widetilde{C}_0\setminus 
\pi^{-1}(\d C_0)$ is an analytic subset 
of $\widehat{D}\setminus \pi^{-1}(\d C_0)$. Let $U_0$ be a relatively compact open subset of $\cc^n$ such that $C_0$ 
is analytic in $U_0$ and $\d C_0= \bar C_0\cap \d U_0$. Let an open $U_1\comp U_0$ be such that $U_0\setminus U_1
\subset W_0$. 
There are only finitely many components of $C_0$ which intersect $\bar U_1$, all other (they can be infinite in 
number) are contained in $U_0\setminus \bar U_1\subset W_0$.
But $\pi|_{\widehat{W}_0}:\widehat{W}_0\to W_0$ is biholomorphic. Therefore there is no problem of lifting these
components to $\widehat{W}_0$. We can forget them and suppose without 
loss of generality that $C_0$ has only finitely many irreducible components, all intersecting $\bar U_1$.

\begin{defi}
\label{lift-c0}
We define $\widehat{C}_0$ to be the union of components of $\widetilde{C}_0$ intersecting $\widehat{W}_0$.
\end{defi}
Recall that by {\slsf components} we mean the irreducible components.
We shall prove that $(\widehat{C}_0, \d\widehat{C}_0)$ satisfies the conclusion of our theorem, where 
$\d\widehat{C}_0$ is defined by \eqqref(1-lift) or, equivalently as $\i (\d C_0)$.

\smallskip\noindent{\slsf Step 3.  Lift of paths to the envelope.} By $\Sing C$ we denote the 
set of singular points of the analytic set $C$ and by $\Reg C = C\setminus \Sing C$ the set of its smooth points.

\begin{lem} 
\label{path-lift0}
Let $\hat p_0, \hat p_1 \in \Reg \widetilde{C}_0$ be points on the same irreducible component of 
$\widetilde{C}_0 \deff \pi^{-1}(C_0)$ such that both $p_0\deff \pi(\hat p_0) \in \Reg C_0$ and 
$p_1\deff \pi(\hat p_1)\in \Reg C_0$. Let $\hat\gamma = \{\hat\gamma (\tau ): \tau \in [0,1]\} $ 
be a continuous path in $\Reg \widetilde{C}_0$ from $\hat p_0$ to $\hat p_1$ such that 
$\gamma (\tau )\deff \pi (\hat\gamma (\tau )) \in \Reg C_0$ for all $\tau \in [0,1]$. Suppose 
that for every $k\gg 1$ there exists $p_k\in \Reg C_k$ such that 

\sli $p_k\to p_0$ as $k\to\infty $;

\slii $\i (p_k) \to \hat p_0$ as $k\to\infty $.

\noindent Then for $k\gg 1$ there exists a path $\gamma_k = \{ \gamma_k(\tau ): \tau \in [0,1]\}$ in 
$\Reg C_k$ such that its canonical lift $\i (\gamma_k)$ to $\widehat{D}$ is arbitrarily close to 
$\hat\gamma$. 
\end{lem}
\proof More precisely we mean that for any $\eps >0$ there exists $k_0$ such that for every $k\ge k_0$ there exists a path $\gamma_k$ in $\Reg C_k$ such that

\begin{equation}
\eqqno(path-dist1)
\dist (\hat\gamma , \i (\gamma_k))\deff \sup\{d(\hat\gamma (\tau ), \i (\gamma_k(\tau ))): \tau \in [0,1]\} < \eps .
\end{equation}
Note that $C_k\ni p_k\to p_0$ always exist simply because $\dist_H(\bar C_k, \bar C_0)\to 0$. In our applications of this lemma we shall always we able to find $p_k$-s such that $\i (p_k)\to \hat p_0$ as
well. 

\smallskip Perturbing our path we can assume that $\hat\gamma $ is real analytic, has only transverse 
self-intersections if $\dim\widetilde{C}_0=1$ or, is imbedded if $\dim\widetilde{C}_0 > 1$. Moreover, 
we can assume  that the same holds for its projection $\gamma $ in $C_0$. Take a tubular Stein neighborhood 
$U$ of $\gamma $ in $C_0$ such that $\bar U \subset \Reg C_0$. In the case $\dim C_0=1$ it is a 
self-intersecting band, see the Picture \ref{band} below. Let $N$ be a Stein neighborhood of $U$ in $\cc^n$ which is biholomorphic to a neighborhood of the zero section in the normal 
bundle to $U$, see \cite{Siu} or \cite{De}. Denote by $\pi_N:N\to U$ the natural projection thus obtained.
If $k_0$ is taken sufficiently big we have for every $k\ge k_0$ that  $C_k\cap \d N \subset \pi_N^{-1}(\d U)$  and therefore $\pi_N|_{C_k}:C_k\cap N \to U$ is proper, consequently is an analytic cover.

\smallskip Appropriately perturbing $\hat\gamma $ together with its $\pi$-projection $\gamma$ we can assume that the path $\gamma$ doesn't passes through the branch locus of $\pi_N|_{C_k}$ and stays to be real analytic.

\begin{rema} \rm
\label{sheet-numb}
The number of sheets of this cover might be non bounded when $k\to \infty $. But for every $k$ the set $B_k$ of 
branch points of $\pi_N|_{C_k}:C_k\cap N\to U$ is a finite subset of $U$. For this one may need to shrink $U$, 
and then this will hold for every $k\gg 1$. Therefore it is clear that every path, say $\gamma $ can be  
approximated by real analytic paths avoiding $B_k$. We need this only for a fixed $k$ in fact.
\end{rema}

\begin{figure}[h]
\centering
\includegraphics[width=3.0in]{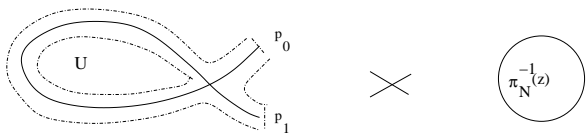}
\caption{After shrinking $U$ we can assume that $\pi_N:N\to U$ is a $(n-q)$-disk bundle over $U$ (not
necessarily trivial as on this picture) and its boundary consists from two obvious parts: first is 
$\pi_N^{-1}(\d U)$ (which is diffeomorphic to a disk bundle over the boundary $\d U$), 
and the second is $\overline{\d N\setminus \pi_N^{-1}(\d U)} = \bigcup_{z\in \bar U}\d \pi_N^{-1}(z)$. 
The latter is a circle (sphere if $(n-q)>1$) bundle over $\bar U$.}
\label{band}
\end{figure}

Using the fact that 
\begin{equation}
\eqqno(lift-k)
\pi_N|_{C_k}: (C_k\cap N)\setminus \pi_N|_{C_k}^{-1}(B_k)\to U\setminus B_k 
\end{equation}
is an unramified cover we can perturb $p_k$ on $C_k$ and $p_0$ on $C_0$ in order that $p_k\in \pi_N^{-1}(p_0)$.
Indeed, we can do this since $p_k$ is close to $p_0$ and therefore $\pi_N(p_k)$ is close to $p_0$ as well.
Now we lift uniquely $\gamma $ to a path $\gamma_k\subset C_k$ under this covering with 
initial point $\gamma_k(0) = p_k$. We are going to prove that the natural 
inclusion $\i (\gamma_k(\tau ))$ is close to $\hat\gamma (\tau )$ for all $\tau \in [0,1]$. Remark that the polydisk distance 
between $\gamma_k(\tau )$ and $\gamma (\tau )$ is not more than $d_1$ for all $\tau \in [0,1]$ if $k$ was taken sufficiently 
big.  Let $V_{\tau}$ be the connected component of $\Delta^n(\gamma (\tau ),d)\cap D$ containing $\gamma_k(\tau )$. Notice 
that by Proposition \ref{behn-pr3} we have that for every $f\in \calo (D)$ the restriction $f|_{V_{\tau}}$ holomorphically 
extends to $\Delta (\gamma (\tau ), d)$.

\begin{claim}
\label{claim-3.1}
We claim that for any $f\in \calo(D)$ extensions 
$\ext (f|_{V_{\tau}})$ and $\ext (f|_{V_{\nu}})$ coincide on $\Delta^n(\gamma (\tau ), d)\cap \Delta^n(\gamma (\nu),d)$ 
provided $|\tau -\nu|$ is small enough.
\end{claim}
Since $\bar C_k\comp D$ we can find some $0<d_2<d_1$ such that $\dist (C_k,\d D) > d_2$. 
Therefore both $\Delta^n(\gamma_k (\tau ), d_2)\subset V_{\tau }$ and $\Delta^n(\gamma_k (\nu),d_2)\subset V_{\tau}$, and as a 
consequence $V_{\tau }\cap V_{\tau}\supset \Delta^n(\gamma_k(t), d_2)\cap \Delta^n(\gamma_k(\nu),d_2)\subset D$. The latter 
intersection is non-empty provided  $\dist (\gamma_k(\tau ), \gamma_k (\nu)) < d_2$, \ie if $|\tau -\nu |$ is small enough. 
Since for every $\tau \in [0,1]$ the extension $\ext (f|_{V_{\tau}})$ was obtained as an extension of $f$ from a neighborhood 
of $\gamma_k (\tau )$ the claim follows from the uniqueness theorem for holomorphic functions.

\smallskip  By a lift of $\Delta^n(\gamma (\tau ),d)$ to $\hat D$ we understand a biholomorphism $l_{\tau }$ of $\Delta^n
(\gamma (\tau ),d)$ onto some domain in $\hat D$  such that $\pi\circ l_{\tau } = \id$. This domain we shall denote as $\hat\Delta^n
(\hat\gamma (\tau ),d)$ and shall justify this notation by proving that in our case $\hat\Delta^n(\hat\gamma (\tau ),d)\ni 
\hat\gamma (\tau )$ and that $\pi|_{\hat\Delta^n(\hat\gamma (\tau ),d)} : \hat\Delta^n(\hat\gamma (\tau ),d)\to \Delta^n(\gamma (\tau),d)$
is a biholomorphism sending $\hat \gamma(\tau )$ to $\gamma (\tau)$ for all $\tau \in [0,1]$. The lifts in question are given by simultaneous 
extensions of all holomorphic in $D$ functions from $V_{\tau }$ to $\Delta^n(\gamma (\tau ),d)$ via Proposition \ref{behn-pr3} of Step 1,
and then taking corresponding homomorphisms. 


\begin{claim}
\label{claim-3.2}
The lifts $l_{\tau }:\Delta^n(\gamma (\tau ),d)\to \hat\Delta^n(\hat\gamma (\tau ),d)$ of polydisks 
$\Delta^n(\gamma (\tau ),d)$ to $\hat D$ constructed as above are coherent in the sense that

\smallskip a) $\hat\Delta^n(\hat\gamma (\tau ),d)\cap \hat\Delta^n(\hat\gamma (\nu),d) \not= \emptyset$ for 
$|\tau -\nu |$ small enough. 

\smallskip b) Moreover, $\hat\Delta^n(\hat\gamma (\tau ),d)$ contains $\hat\gamma (\tau )$ for every 
$\tau \in [0,1]$.
\end{claim}
Notice that for every $\tau \in [0,1]$ the biholomorphic lift of $\Delta^n (\gamma_k(\tau ),d + d_1)$ to 
$\widehat{D}$ coincides on $\Delta^n (\gamma_k(\tau ),d_2)$ with the canonical lift $\i : \Delta^n (\gamma_k(\tau ),d_2) 
\to \hat\Delta^n (\i(\gamma_k(\tau )),d_2)$. The first part of the  claim is obvious since extensions of $f\in \calo (D)$ 
to $\Delta^n(\gamma (\tau ),d)$ are extensions from $\Delta^n (\gamma_k(\tau ),d_2)$ to $\Delta^n (\gamma_k(\tau ),d + d_1)
\supset \Delta^n(\gamma (\tau ), d)$, where $d_2$ was taken small enough, namely such that $\Delta^n(\gamma_k(\tau ),d_2)\subset D$. 
Second, since the distance between $\i (p_k)$ and $\hat p_0$ is supposed to be small we get that 
$\hat\Delta^n(\i (p_k),d)\ni \hat p_0$. Therefore $\hat\Delta^n(\i (\gamma_k(0)),d)\ni \hat \gamma (0)$ and therefore 
$l_0(\gamma (0)) = \hat\gamma (0)$.  What we need to prove is that $l_{\tau }(\gamma (\tau )) =\hat\gamma (\tau )$ for 
all $\tau \in [0,1]$. For the starting lifts $l_{\tau } : \Delta^n(\gamma (\tau ),d) \to \hat\Delta^n(\hat\gamma (\tau ),
d)$, $\tau \sim 0$, this is true by local biholomorphicity of $\pi$. And therefore $\hat\Delta^n(\hat\gamma (\tau ),d)
\ni \hat\gamma (\tau )$, \ie $l_{\tau }(\gamma (\tau )) = \hat\gamma (\tau )$ for $\tau \sim 0$. By real analyticity of
both $\gamma $ and $\hat\gamma $ the relation $l_{\tau }(\gamma (\tau )) = \hat\gamma (\tau )$ stays valid for all 
$\tau \in [0,1]$. The claim is proved.

\smallskip By construction the polydisk distance between $\i (\gamma_k (\tau ))$ and $\hat\gamma (\tau )$  is not more 
than $d_1$ because this is true for $\gamma_k (\tau )$ and $\gamma (\tau )$. Taking  $0< d_1 <\eps$ and $k_0$ as above
we get the proof of the lemma.

\smallskip\qed

In a particular case when we consider an irreducible component of $\widehat{C}_0$ and when the starting point $\hat p_0 
\in \widehat{W_0}\cap \Reg \widehat{C}_0$, \ie $\hat p_0$ is close to the boundary, such $p_k$ obviously exists because
$\d C_k\to \d C_0$ in Hausdorff topology. Moreover, since $\i$ locally preserves the polydisk distance, we have that 
$\i (p_k)\to \hat p_0$ as well. Therefore we obtain the following

\begin{corol}
\label{path-lift1}
Let $\hat p_0 \in \widehat{W_0}\cap \Reg \widehat{C}_0$ and $\hat p_1 \in \Reg \widehat{C}_0$ be points on the same 
irreducible component of $\widehat{C}_0$ such that $p_0\deff \pi(\hat p_0) \in W_0\cap \Reg C_0$ and $p_1\deff 
\pi(\hat p_1)\in \Reg C_0$. Let $\hat\gamma \deff \{\hat\gamma (\tau ): \tau \in [0,1]\} $ be a continuous path in 
$\Reg \widehat{C}_0$ from $\hat p_0$ to $\hat p_1$ such that $\gamma (\tau )\deff \pi (\hat\gamma (\tau )) \in 
\Reg C_0$ for all $\tau \in [0,1]$. Then for $k\gg 1$ there exists a path $\gamma_k\deff \{ \gamma_k(\tau ): \tau 
\in [0,1]\}$ in $\Reg C_k$ such that its canonical lift $\i (\gamma_k)$  to $\widehat{D}$ is arbitrary close to 
$\hat\gamma$. 
\end{corol}

\smallskip\noindent{\slsf Step 4. Projection $\pi|_{\widehat{C}_0}: \widehat{C}_0\to C_0$ is proper.}
This will follow from the following statement.

\begin{lem}
For any point $\hat p_1\in \widehat{C}_0$ and any $f\in \calo (D)$ one has 
\begin{equation}
\eqqno(max-pr2)
|\hat f(\hat p_1)| \le 2^nM(f, \overline{\d C_0^{d+2d_1}}).
\end{equation}
Here $\hat f$ is the canonical extension of $f$ to the envelope $\widehat{D}$.
\end{lem}
\proof Perturbing $\hat p_1$ a bit we can assume that $\hat p_1 \in \Reg\widehat{C}_0$, as well as $p_1=\pi (\hat p_1) 
\in \Reg C_0$.  Now take $\hat p_0$ on the same irreducible component of $\hat C_0$ as $\hat p_1$ and such that 
$\hat p_0$ is close to $\d\widehat{C}_0$. In addition take $\hat p_0$ such that both $\hat p_0$ and $p_0=\pi (\hat p_0)$ 
lie on the smooth locuses of $\widehat{C}_0$ and $C_0$ respectively. Take a path $\hat\gamma $ from  $\hat p_0$ to 
$\hat p_1$ in $\Reg\widehat{C}_0$ and a path $\gamma_k$ in $C_k$ as in Corollary \ref{path-lift1}. From the fact that 
$\hat\gamma (\tau )$ belongs to $\hat\Delta^n(\i (\gamma_k(\tau )),d_1)$ for every $\tau $ we get from \eqqref(cauchy2)
the following estimate
\begin{equation}
\eqqno(cauchy3)
|\hat f(\hat\gamma (\tau ))| \le  \sum_{|m|=0}^{\infty}|a_m(\hat f, \i(\gamma_k(\tau )))| \cdot ||\hat\gamma (\tau ) - 
\i (\gamma_k(\tau ))||^{|m|} \le \sum_{|m|=0}^{\infty}|a_m(f,\gamma_k(\tau ))|d_1^{|m|}\le
\end{equation}
\[
\le M(f, \overline{\d C_0^{d+2d_1}})
\sum_{|m|=0}^{\infty}\frac{d_1^{|m|}}{(d+d_1)^{|m|}} = M(f, \overline{\d C_0^{d+2d_1}})\left(1 + \frac{d_1}{d}\right)^n.
\]
This implies \eqqref(max-pr2) and lemma is proved.

\smallskip\qed 

As for the properness of $\pi|_{\widehat{C}_0}: \widehat{C}_0\to C_0$ proceed as follows. Since, as it was explained 
above, we can suppose that $C_0$ has only finitely many irreducible components (say it is itself irreducible), it is 
sufficient to prove that $\pi|_{\widehat{C}_0} : \widehat{C}_0\to C_0$ is proper for an irreducible $\widehat{C}_0$. 
If not there would exist a sequence of points $p_j\in C_0$ converging to a point $p_0\in C_0$ and a sequence 
of their $\pi$-preimages $\hat p_j\in \widehat{C}_0$ leaving every compact in $\widehat{C}_0$. Would $\hat p_j$ leave every 
compact in $\widehat{D}$ this would contradict to \eqqref(max-pr2) and the holomorphic convexity of $\widehat{D}$. 
Otherwise, modulo taking a subsequence, we would have that $\hat p_j$ converge to a point $\hat p_0\in \pi^{-1}(p_0)$ 
which is not in $\widehat{C}_0$. But this is not possible since $\pi$ is biholomorphic between neighborhoods of 
$\hat p_0$ and $p_0$. 

\begin{rema} \rm
\label{bihol-pi}
Notice that \eqqref(max-pr2) applied to the Taylor coefficients $a_m(f,p)$ of holomorphic in $D$ functions implies 
that 
\[
\widehat{C}_0\subset \widehat{D}_{d_0} \deff \{p\in \widehat{D}: \dist (p, \d\widehat{D}) > d_0\},
\]
and, moreover, that for every $p\in C_0$ and every $\hat p\in \pi^{-1}(p)\cap \widehat{C}_0$ the restriction 
\[
\pi|_{\hat\Delta^n(\hat p,d)} : \hat\Delta^n(\hat p,d) \to \Delta^n(p, d)
\]
is a biholomorphism. Here $0<d<d_0$ can be taken as close to $d_0$ as one wishes, but the component $V$ of 
$\Delta^n(p, d)\cap D$ from which all functions $f\in \calo(D)$ extend to $\Delta^n(p,d)$ may not be the same for 
all $d$. This follows from the Claim \ref{claim-3.2} with $\hat p = \hat\gamma (1)$ and 
$p=\gamma (1)$ for an appropriate paths $\hat\gamma$ on $\hat C_0$ and $\gamma = \pi (\hat\gamma)$ on $C_O$. 
The eventual dependence of $V$ from $d$ appears when we choose the approximating path $\gamma_k$ on $C_k$.
\end{rema}

\medskip\noindent{\slsf Step 5. $\widehat{C}_0$ is an analytic set in $\widehat{D}\setminus \d\widehat{C}_0$ and its 
boundary is $\i (\d C_0)$.}
Up to now we know that $\widehat{C}_0$ is an analytic set in $\widehat{D}\setminus \pi^{-1}(\d C_0)$. What
we need to prove is that $\widehat{C}_0$ cannot accumulate to $\pi^{-1}(\d C_0)\setminus \i (\d C_0)$. Suppose 
this is wrong. Then one can find $\hat p_1$ close to $\pi^{-1}(\d C_0)\setminus \i (\d C_0)$ and $\hat p_0$ near 
$\d\widehat{C}_0 = \i (\d C_0)$ which belongs to the same irreducible component of $\widehat{C}_0$ as $\hat p_1$,
both projecting to points $p_1, p_0\in C_0$, where $p_0$ is near $\d C_0$. But notice that from the properness of 
$\pi|_{\widehat{C}_0}: \widehat{C}_0\to C_0$ it follows that $p_1$ is also near $\d C_0$. After the obvious 
perturbations take a path $\hat\gamma $ from $\hat p_0$ to $\hat p_1$ on the regular part of $\widehat{C}_0$. Let 
$\gamma_k$ be a path on $C_k$ as in Corollary \ref{path-lift1} \ie its canonical lift  $\i (\gamma_k)$ is close to 
$\hat\gamma $. Since $\gamma_k(0) \sim p_0$ and $\gamma_k(1) \sim p_1$ and they are both close to $\d C_0$ their canonical 
lifts are both close to $\i (C_0)=\d\widehat{C}_0$. Therefore both $\hat\gamma (0)$ and $\hat\gamma (1)$ are close to 
$\d \widehat{C}_0$. For the case of $\hat\gamma (1)$ this is a contradiction.

\smallskip\noindent{\slsf Step 6. $\i (C_k)$ converge to $\widehat{C}_0$.}  In addition to the Remark \ref{bihol-pi}
above we see that $\widehat{C}_0$ is contained in a bounded part of 
$\widehat{D}$, this follows from the properness of $\pi|_{\widehat{C}_0}: \widehat{C}_0 \to C_0$. Therefore 
$(\widehat{C}_0, \d\widehat{C}_0)$ is a compact analytic set with boundary in an appropriately taken $\widehat{U}
\comp \widehat{D}$. Moreover, for every $k$ we have that $\i (C_k)$ is a compact analytic set with boundary in 
$\i (U_k)\subset \hat D$ for an appropriate neighborhood $U_k$ of $C_k$.

\smallskip Let $\hat p_1$ be any point of $\widehat{C}_0$. We need to approximate it by points $\i (p_k)$ with 
$p_k\in C_k$. Perturbing $\hat p_1$ slightly we can suppose that $\hat p_1\in \Reg \widehat{C}_0$. Fix some 
$\hat p_0\in \Reg\widehat{C}_0 \cap \widehat{W}_0$ on the same irreducible component of $\widehat{C}_0$ as 
$\hat p_1$. Fix some path $\hat\gamma$ from $\hat p_0$ to $\hat p_1$ on $\Reg\widehat{C}_0$. By Corollary 
\ref{path-lift1} we can approximate $\hat\gamma$ by $\i (\gamma_k)$ with $\gamma_k\subset \Reg C_k$. Now 
$\i (\gamma_k(1))$ will approximate $\hat p_1$. Therefore $\lim_H\i (C_k)\supset \widehat{C}_0$. 

\smallskip To prove the opposite inclusion suppose that the Hausdorff limit of $\i (C_k)$ is bigger than 
$\widehat{C}_0$, \ie that there exists $\hat p_0\in \widetilde{C}_0\setminus \widehat{C}_0$ which is in 
$\lim_H\i (C_k)$ lying on some irreducible component $C'$ of $\widetilde{C}_0$ which is not in $\widehat{C}_0$.
There exist $C_k\ni p_k\to p_0\deff \pi (\hat p_0)$ such that $\i (p_k)\to \hat p_0$. Take $\hat p_1$ on the 
same irreducible component of $\widetilde{C}_0$ as $\hat p_0$ close to its boundary $\d\tilde C_0$. Perturbing slightly 
all points in question we assume that they lie on the smooth locuses of corresponding analytic sets. Take a path
$\hat\gamma \subset \Reg\widetilde{C}_0$ from $\hat p_0$ to $\hat p_1$ and approximate it by $\i (\gamma_k)$ for 
$\gamma_k\subset \Reg C_k$ as in Lemma \ref{path-lift0}. If $\hat p_1$ was close to $\pi^{-1}(\d C_0)$ than 
$\gamma_k(1)$ must be close to $\d C_0$ and therefore $\i (\gamma_k(1))$ will be close to $\i (\d C_0) = 
\d\widehat{C}_0$. This implies in its turn that $\hat p_1$ is close to $\i (\d C_0)$ contradicting to the assumption 
that the irreducible component $C'$ we working with doesn't belong to $\widehat{C}_0$. Therefore $\hat p_1$ must 
be close to infinity in $\hat D$. But then by holomorphic convexity of $\hat D$ we find a holomorphic function $f$ 
on $D$ such that $|\hat f(\hat p_1)| > \{\sup |f(p)|:p\in \d C_0\}$. The same will hold true for $\hat f(\gamma_k(1))$ 
for $k$ big enough. This contradicts to the maximum principle for holomorphic functions on $C_k$ or, equivalently on 
$\i (C_k)$. The step is proved.

\smallskip Theorem is proved in the case $X=\cc^n$.

\smallskip\noindent{\slsf Step 7. Case of Stein $X$.} Now consider the case when $X$ is a Stein manifold.
Imbed $X$ to $\cc^N$ properly and let $\ret :V\to X$ be a holomorphic retraction of an appropriate 
neighborhood of $X$. Set $D^{\ret}\deff \ret^{-1}(D)$. Then $C_k$-s and $C_0$ are clearly compact analytic sets
with boundary in $D^{\ret}$. We can repeat the consideration as above for the algebra $\calo^{\ret}(D^{\ret})$ of 
$\ret$-invariant holomorphic functions on $D^{\ret}$ and get the lift $\widehat{C}^{\ret}_0$ to the $\ret$-invariant 
envelope $(\widehat{D}^{\ret}, \pi^{\ret})$. Obviously $\widehat{C}^{\ret}_0\subset (\pi^{\ret})^{-1}(X)$ and satisfies 
the conclusions of the theorem. 

\smallskip\qed

\newprg[LIFT.eq-d]{Proof of Proposition \ref{cp1}.}

We need to prove that if $q \ge n/2$  the lift $\widehat{C}_0$ of  $C_0$ is singlevalued. First let us 
prove the following property of the projection $\pi$ near $\widehat{C}_0$. Take some $0<d<d_0 = \dist (\d C_0, \d D)$. 
We know that for every $p\in C_0$ and every $\hat p\in \pi^{-1}(p)\cap \widehat{C}_0$ the restriction 
$\pi|_{\hat\Delta^n(\hat p,d)} : \hat\Delta^n(\hat p,d) \to \Delta^n(p, d)$ is a biholomorphism, see Remark \ref{bihol-pi}.

\begin{lem}
\label{lift-uniq}
For every $p\in C_0$ and every irreducible component $C$ of $C_0\cap \Delta^n(p,d)$ there exists a unique $\hat p\in 
\pi^{-1}(p)$ such that $C$ lifts biholomorphically to an irreducible component $\widehat{C}$ of $\widehat{C}_0\cap 
\hat\Delta^n(\hat p , d)$. 
\end{lem}
\proof Suppose that this is wrong. Then there exist two distinct points $\hat p_1$ and $\hat p_2$ in $\widehat{C}_0$ 
such that $\pi (\hat p_1) = \pi (\hat p_2)=p\in C_0$ and components $\widehat{C}^{'}$ of $\widehat{C}_0\cap 
\hat\Delta^n(\hat p_1 , d)$ and $\widehat{C}^{''}$ of $\widehat{C}_0\cap \hat\Delta^n(\hat p_2 , d)$ respectively 
which are mapped by $\pi$ onto the same component $C$ of $C_0\cap \Delta^n(p,d)$. Perturbing $p$ slightly we can 
assume that $p\in \Reg C_0$ as well as $\hat p_1, \hat p_2\in \Reg\hat C_0$. Take a smooth path $\gamma $ in 
$\Reg C_0$ from $p$ to some $q\in \Reg C_0\cap W_0$. Using the fact that $\pi|_{\widehat{C}_0}: \widehat{C}_0
\setminus \pi^{-1}(\Sing C_0) \to \Reg C_0$ is an analytic cover we can lift $\gamma$ to $\hat\gamma_1\subset 
\Reg\widehat{C}_0$ starting at $\hat p_1$ and to $\hat\gamma_2\subset \Reg\widehat{C}_0$ starting at $\hat p_2$. 
These paths end at $\hat q_1\in \widehat{C}_0\cap\widehat{W}_0$ and $\hat q_2\in \widehat{C}_0\cap\widehat{W}_0$ respectively. But  our projection $\pi$ is biholomorphic when restricted to $\widehat{C}_0\cap \widehat{W}_0$. 
Therefore $\hat q_1=\hat q_2 = \i (q)$. And this implies that $\hat p_1 = \hat p_2$. Contradiction. 
Lemma is proved.

\smallskip\qed

\smallskip Now suppose that the lift $\widehat{C}_0$ of our $C_0$ is not singlevalued, \ie that there exist two 
distinct points $\hat p_1$ and $\hat p_2$ in $\widehat{C}_0$ such that $\pi (\hat p_1) = \pi (\hat p_2)=p\in C_0$. 
From the lemma just proved it follows that $p$ belongs to the intersection of two {\slsf distinct} irreducible 
components $C'$ and $C^{''}$ of $C_0\cap \Delta^n(p,d)$ and there exist irreducible components $\widehat{C}^{'}$ 
of $\widehat{C}_0\cap \hat\Delta^n(\hat p_1,d)$ and $\widehat{C}^{''}$ of $\widehat{C}_0\cap \hat\Delta^n(\hat p_2,d)$ 
which are mapped by $\pi$ onto $C'$ and $C{''}$ respectively.

\begin{claim} 
\label{claim-3.3}
One can find irreducible components $\widehat{C}_k'$  of $\widehat{C}_k\cap \hat\Delta^n(\hat p_1,d)$ and 
$\widehat{C}_k^{''}$ of $\widehat{C}_k\cap \hat\Delta^n(\hat p_2,d)$   such that the Hausdorff limit of 
$\widehat{C}_k'$ contains $\widehat{C}'$ and the Hausdorff limit of $\widehat{C}_k^{''}$ contains $\widehat{C}'$. 
\end{claim}

Indeed, the Hausdorff limit of $\widehat{C}_k\cap \hat\Delta^n(\hat p_i,d)$ contains $\widehat{C}_0\cap 
\hat\Delta^n(\hat p_i,d)$ for $i=1,2$. Take two points $p'\in \Reg C'$ and $p^{''}\in \Reg C^{''}$. In an
appropriate local coordinates near $p'$ represent $C'$ as $z'=0$, where $z'=(z^{q+1}, ..., z^n)$.
For $k\gg 1$ the piece of $C_k$ will be an analytic  cover of $\Delta^q$ (everything in these local 
coordinates). The claim follows.

\smallskip Therefore, adding to $\widehat{C}^{'}$ and $\widehat{C}^{''}$ some more components we have that 
\[
\lim_H\widehat{C}_k' = \widehat{C}^{'} \qquad \text{ and }\qquad \lim_H\widehat{C}_k' = \widehat{C}^{''}.
\]

\smallskip If we suppose now that $C_k'\deff \pi (\widehat{C}_k')$ and $C_k^{''}\deff \pi (\widehat{C}_k^{''})$ 
are equal for an infinite number of $k$-s (or, even for some $k$-s) then it leads to a contradiction with the 
fact that $\i (C_k)$ is imbedded to $\widehat{D}$. Therefore $C_k'$ and $C_k^{''}$ are distinct for $k\gg 1$. 
Since their Hausdorff limit is $C'$ and $C^{''}$ respectively, which intersect at least at $p$, and due to 
the assumption that $\dim C_k\ge n/2$ we conclude that $C_k'$ intersect $C_k^{''}$ for $k\gg 1$. Indeed, we 
can consider two cases.

\smallskip\noindent{\slsf Case 1.} {\it $\dim C'\cap C^{''} = 0$, in particular $p$ is an isolated point of this intersection.}
Take $d >0$ small enough in order that $p$ is the only point in $C'\cap C^{''}$ and consider the analytic set 
$A\deff C'\times C^{''}$ in $\Delta^n(p,d)\times \Delta^n(p,d)$. It intersects diagonal $\dd$ of $\Delta^n(p,d)
\times \Delta^n(p,d)$ by exactly one point, namely by $p^2\deff (p,p)$. Therefore the restriction to $A$ of the 
projection $p:\Delta^n(p,d)\times \Delta^n(p,d)\to \dd^{\perp}$ to the orthogonal $\dd^{\perp}$ parallel to $\dd$ 
is proper, \ie is an analytic cover. Notice that since $C_k^{'}$ (resp. $C_k^{''}$) converges to $C'$ (resp. to 
$C^{''}$) we have that $A_k\deff C_k'\times C_k^{''}$ converges to $A$. Therefore for $k\gg 1$ $p|_{A_k}:A_k\to 
\dd^{\perp}$ is an analytic cover as well and therefore intersects the vertical $\{p\}\times \dd$. I.e., intersects 
the diagonal, say by $p_k^2\deff (p_k,p_k)$. But then this $p_k$ is a point of intersection of $C_k^{'}$ with 
$C_k^{''}$.

\smallskip\noindent{\slsf Case 2.} {\it $\dim C'\cap C^{''} >0$.} The analytic set $C'\cap C^{''}$ having positive 
dimension, reaches the boundary of $C_0$, \ie intersects $W_0$. At a generic point on this analytic set two local 
branches of $C_0$ should intersect. Moreover, they are biholomorphic images of two nonintersecting local branches 
of $\widehat{C}_0$. Near $\d C_0$ this  contradicts to the fact that $C_0$ lifts biholomorphically to $\i (C_0)$ 
near the boundary. Proposition \ref{cpd} is proved.

\smallskip\qed

\newprg[LIFT.cp2]{Proof of Proposition \ref{cp2}}

Set 

\begin{equation}
\eqqno(d-nol-1)
d_0 = \dist (\bigcup_{t\in [0,1]}\d C_t, \d D). 
\end{equation}
Remark that the conditions of the theorem imply that $d_0>0$. Take a neighborhood $W_0$ of 
$\bigcup_{t\in [0,1]}\d C_t$ such that for $\widehat{W}_0 = \i (W_0)$ the restriction 
$\pi|_{\widehat{W}_0}:\widehat{W}_0\to W_0$ is a biholomorphism. It will be convenient for us to 
prove together with Proposition \ref{cp2} also the following additional 

\begin{asser}
\label{asser1}
Set $\d\widehat{C}_t\deff \i (\d C_t)$ for every $t\in [0,1]$ and let $\widehat{C}_{t}\subset 
\widehat{D}$ be the union of all irreducible components of $\widetilde{C}_{t} \deff \pi^{-1}(C_{t})$ 
that intersect $\widehat{W}_0$. Then $\{(\widehat{C}_t, 
\d\widehat{C}_t)\}_{t\in [0,1]}$ is the family satisfying the conclusion of Proposition \ref{cp2}. In
particular $\i (\d C_t)$ is the whole boundary of $\widehat{C}_t$, justifying the notation $\d \widehat{C}_t
= \i (\d C_t)$ above.
\end{asser}

\medskip
Denote by $T$ the set of $t \in [0,1]$ such $\{ (C_t, \d C_t)\}_{[0,1]}$ can be {\slsf continuously} lifted to 
$\hat D$ up to $t$. By saying that we mean that the family $\{(\widehat{C}_t,\d\widehat{C}_t)\}_{[0,t]}$ 
constructed as in Assertion \ref{asser1} satisfies the conclusion of Theorem \ref{cp2}. $T$ is non-empty 
since it contains a neighborhood of zero.

\smallskip\noindent{\slsf $T$ is closed.} Let $t_0=\sup \{t:t\in T\}$. We need to prove that $t_0\in T$.
The proof of this statement follows the main lines of the proof of Theorem \ref{cpd}.
The following lemma is analogous to Lemma \ref{path-lift0}. As in the quoted lemma we 
can assume that the number of irreducible components of $C_{t_0}$  is finite.

\begin{lem}
\label{path-lift2}
Let $\hat p_0 , \hat p_1 \in \Reg \widehat{C}_{t_0}$ be points on the same irreducible component of $\widehat{C}_{t_0}$ 
such that $p_0\deff \pi(\hat p_0) \in  \Reg C_{t_0}$ and $p_1\deff \pi(\hat p_1)\in \Reg C_{t_0}$. Let  $\hat\gamma = 
\{\hat\gamma (\tau ): \tau\in [0,1]\} $ be a path on $\Reg \widehat{C}_{t_0}$ from $\hat p_0$ to $\hat p_1$ such that 
$\gamma (\tau )\deff \pi (\hat\gamma(\tau )) \in \Reg C_{t_0}$ for all $\tau \in [0,1]$. Suppose that for $t\sim t_0$

\sli there exist $p_t\in \Reg C_t$ such that $p_t\to p_0$ as $t\to t_0$;

\slii there exist lifts $\hat p_t\in \Reg\widehat{C}_t$ of $p_t$ such that $\hat p_t\to \hat p_0$ as $t\to t_0$.

\noindent Then for $t\sim t_0$ there exists a path $\gamma_t = \{ \gamma_t(\tau ):\tau\in [0,1]\}$ 
on $\Reg C_t$ such that its lift $\hat\gamma_t\subset \Reg\widetilde{C}_t \subset \widehat{D}$ is close to $\hat\gamma$. 
\end{lem}

Since in the case when $\hat p _0\in \widehat{W}_0\cap \Reg\widehat{C}_{t_0}$ the existence of $p_t\in \Reg C_t$ for 
$t\sim t_0$ satisfying items (\sli and (\slii is automatic we have, as in the case of Lemma \ref{path-lift0}, the following
corollary. Here as $\hat p_t$ one should take $\i (p_t)\in \Reg\widehat{C}_t$. One remarks as well that the lifts 
$\hat \gamma_t$ are contained in $\widehat{C}_t$ in this case.

\begin{corol}
\label{path-lift3}
Let $\hat p_0 \in \widehat{W}_0\cap \Reg \widehat{C}_{t_0}$ and $\hat p_1 \in \Reg \widehat{C}_{t_0}$ be points 
on the same irreducible component of $\widehat{C}_{t_0}$ such that $p_0\deff \pi(\hat p_0) \in W_0\cap \Reg C_{t_0}$ and 
$p_1\deff \pi(\hat p_1)\in \Reg C_{t_0}$. Let furthermore $\hat\gamma = \{\hat\gamma (\tau ): \tau\in [0,1]\} $ be a path on 
$\Reg \widehat{C}_{t_0}$ from $\hat p_0$ to $\hat p_1$ such that $\gamma (\tau )\deff \pi (\hat\gamma(\tau )) \in \Reg 
C_{t_0}$ for all $\tau \in [0,1]$. Then for $t\sim t_0$ there exists a path $\gamma_t = \{ \gamma_t(\tau ):\tau\in [0,1]\}$ 
on $\Reg C_t$ such that its lift $\hat\gamma_t\subset \Reg\widehat{C}_t \subset \hat D$ is close to $\hat\gamma$. 
\end{corol}

\proof {\slsf Case 1.} {\it First we shall prove this lemma for $t\lesssim t_0$.} Here when writing 
$t\lesssim t_0$ we mean $t<t_0$ and $t\sim t_0$ \ie close to and less than $t_0$. 
Take some $0 < d < d_0$ close to $d_0$ and some $0< d_1 < \min\{\frac{d_0-d}{4},d\}$. Take $t\lesssim t_0$ such that 
\[
\dist_H (C_{t}, C_{t_0}) < d_1.
\]
For a holomorphic in $D$ function $f$ denote by $\hat f$ its holomorphic extension to $\hat D$. Remark now that
in the same manner as in \eqqref(cauchy1) and \eqqref(cauchy2) we have for any point $p_t\in C_{t}$ and all $m\in \nn^n$ the 
following estimate
\begin{equation}
\eqqno(cauchy4)
|a_m(\hat f, \hat p_t)| \le \frac{M(f, \overline{\d C_{t_0}^{d+2d_1}})}{(d+d_1)^{|m|}}.
\end{equation}
Here $\hat p_t$ is any lift of $p_t$ to $\widehat{C}_t$. 

\begin{rema} \rm
\label{max-pr1}
Let us underline that \eqqref(cauchy4) holds true due to the stated in the formulation of the proposition (and therefore
assumed here for $t\lesssim t_0$) fact that the boundary of $\widehat{C}_t$ is $\i (\d C_t)$. This is not yet proved
for $t>t_0$, \ie we don't know whether $\i (\d C_t)$ is the whole boundary of $\widehat{C}_t$ when $t>t_0$. Therefore
\eqqref(cauchy4) is established up to now for $t\lesssim t_0$ only.
\end{rema}

Let a neighborhood $U$ of (a real analytically perturbed) path $\gamma$ in $C_{t_0}$ and  $\pi_N : N\to U$ be as in the 
proof of Lemma \ref{path-lift1}. For $t\lesssim t_0$ the restriction $\pi_N|_{C_t\cap N} : C_t\cap N \to U$ is an analytic cover. 
Take as $\gamma_t$ any lift of the path $\gamma $ under $\pi_N|_{C_t\cap N}: (C_t\cap N) \setminus \pi_N|_{C_t\cap N}^{-1}(B_t)
\to U\setminus B_t$, where $B_t$ is the branch locus of this cover. Remark that for $t\lesssim t_0$ 
we have
\begin{equation}
\eqqno(path-dist2)
\dist (\gamma (\tau ), \gamma_t(\tau )) < d_1 \quad\text{ for all } \tau \in [0,1]
\end{equation}
by construction. From \eqqref(cauchy4) we see that for $t\lesssim t_0$ and every $\tau $ the Taylor expansion of $\hat f$ at 
$\widehat{\gamma_t(\tau )}$ has the radius of convergence at least $d+d_1$ for any $f\in \calo (D)$. Here $\widehat{\gamma_t(\tau )}$
is any point in $\widehat{C}_t$ over $\gamma_t(\tau)$. This means that $\pi$ is a biholomorphism between an appropriate polydisks
$\hat\Delta^n(\widehat{\gamma_t (\tau )},d+d_1)\subset \hat D$ and $\Delta^n(\gamma_t (\tau),d+d_1)\subset \cc^n$. This allows us
to construct the coherent lift $\hat\gamma_t $ of the path $\gamma_t$ under $\pi|_{\widehat{C}_t}:\widehat{C}_t\to C_t$ starting
with $\hat\gamma_t(0) =\hat p_t$ for $\gamma_t(0)= p_t$. The distance between $\hat\gamma (\tau )$ and $\hat\gamma_t (\tau)$ is 
not more than $d_1$ by construction. The Case 1 of the lemma is proved. 

\smallskip\qed

\begin{rema} \rm
\label{max-pr2} Estimate \eqqref(cauchy4) by continuity stays valid also for $a_m(\hat f, \hat p_{t_0})$, where $\hat p_{t_0}$ 
is any point in $\widehat{C}_{t_0}$. Indeed, due to the definition of $\widehat{C}_{t_0}$ we can join $\hat p_{t_0}$ with some 
$\hat p_0\in \widehat{W}_0\cap \widehat{C}_{t_0}$ by a path $\hat\gamma$. Approximating this path by paths $\hat \gamma_t\subset 
\Reg\widehat{C}_t$ as in Corollary \ref{path-lift3} (for $t\lesssim t_0$ this corollary is already proved), we obtain that 
$\hat p_{t_0}$ is an accumulation point of $\hat p_t(1)\in \widehat{C}_t$. Therefore \eqqref(cauchy4) follows for $\hat p_{t_0}$ 
from the same estimate for $\hat p_t(1)$ by continuity. This implies in its turn that $\widehat{D}$ contains a $d_0$-neighborhood 
of $\widehat{C}_{t_0}$ in the polydisk norm. 
\end{rema}

\smallskip\noindent{\slsf Case 2.} {\it Now we shall prove the general case. } Notice that we have that $\pi_N|_{C_t\cap N}:C_t\cap N \to U$ 
is an analytic cover for $t\sim t_0$ (not only for $t \lesssim t_0$). Due to Remark \ref{max-pr2} we still have that 
$\hat\Delta^n(\hat\gamma (\tau), d)$ is an imbedded polydisk in $\widehat{D}$ such that 
$\pi|_{\hat\Delta (\hat\gamma (\tau), d)}:
\hat\Delta^n(\hat\gamma (\tau) , d) \to \Delta^n(\gamma (\tau ),d)$ is a biholomorphism for every $0<d<d_0$ and very $\tau 
\in [0,1]$. Let $\gamma_t$ be any lift of (a perturbed)  $\gamma$ under $\pi_N|_{C_t\cap N}: (C_t\cap N) \setminus \pi_N|_{C_t\cap N}^{-1}(B_t)
\to U\setminus B_t$, where $B_t$ is the branch locus of this cover. 
Since $\dist (\gamma (\tau), \gamma_t(\tau)) <d_1$ for every $\tau$ we can coherently lift $\gamma_t$ to the envelope of 
holomorphy using biholomorphisms $\pi|_{\hat\Delta (\hat\gamma (\tau), d)}:\hat\Delta^n(\hat\gamma (\tau) , d) \to \Delta^n(\gamma (\tau ),d)$.
The distance between $\hat\gamma (\tau)$ and $\hat\gamma_t(\tau)$ will be the same as $\dist (\gamma (\tau), \gamma_t(\tau)$, \ie
small. Lemma \ref{path-lift2} is proved.

\smallskip Now one can conclude exactly as in Steps 4 and 5 of the proof of Lemma \ref{path-lift1} that projection 
$\pi_{\widehat{C}_{t_0}}:\widehat{C}_{t_0} \to C_{t_0}$ is proper and that $\widehat{C}_{t_0}$ is an analytic set with 
boundary $\d \widehat{C}_{t_0} \deff \i (\d C_{t_0})$ in $\widehat{D}$. Convergence of $\widehat{C_t}$ to $\widehat{C}_{t_0}$ 
as $t\nearrow t_0$ can be proved exactly as in Theorem \ref{cpd}. The same as in Proposition \ref{cp1} one proves that 
$\pi|_{\widehat{C}_{t_0}}: \widehat{C}_{t_0} \to C_{t_0}$ is an isomorphism since it was the case for $t<t_0$. The closeness 
of $T$ is proved.

\begin{rema} \rm 
\label{non-halph}
Let us notice that up to now we did not use the condition $q\ge \dim X/2$. I.e., for {\slsf any} continuous in Hausdorff 
topology family $\{(C_t,\d C_t)$ of compact analytic sets with boundary as in Theorem \ref{cpd} the set of $t\in [0,1]$
up to which it can be lifted to the envelope is {\slsf closed}! This will be used later in the proof of Theorem \ref{cpc}.
\end{rema}

\smallskip\noindent{\slsf $T$ is open.} It is at this place we need to deploy the assumption that $q\ge n/2$. Example \ref{moy-exmp}
in the next section shows that otherwise this is not true. Since $\pi|_{\widehat{C}_{t_0}}: \widehat{C}_{t_0} \to C_{t_0}$ is an 
isomorphism and $\pi$ is locally biholomorphic we conclude that $\pi$ is a biholomorphism between a neighborhoods of 
$\pi|_{\widehat{C}_{t_0}}: \widehat{C}_{t_0}$ and $C_{t_0}$. And now the claim follows. Proposition \ref{cp2} is proved.

\smallskip\qed

\newsect[LIFT-C]{Lift to the envelope of holomorphy II: continuous case}

\newprg[LIFT-C.exmp]{Example}

\noindent{\slsf Part I.} The first part of our construction will provide us an example to Theorem \ref{cpd}. 
This part is taken from \cite{Iv} and it illustrates that the lift in Theorem \ref{cpd} is not singlevalued 
in general.

\begin{figure}[h]
\centering
\includegraphics[width=2.2in]{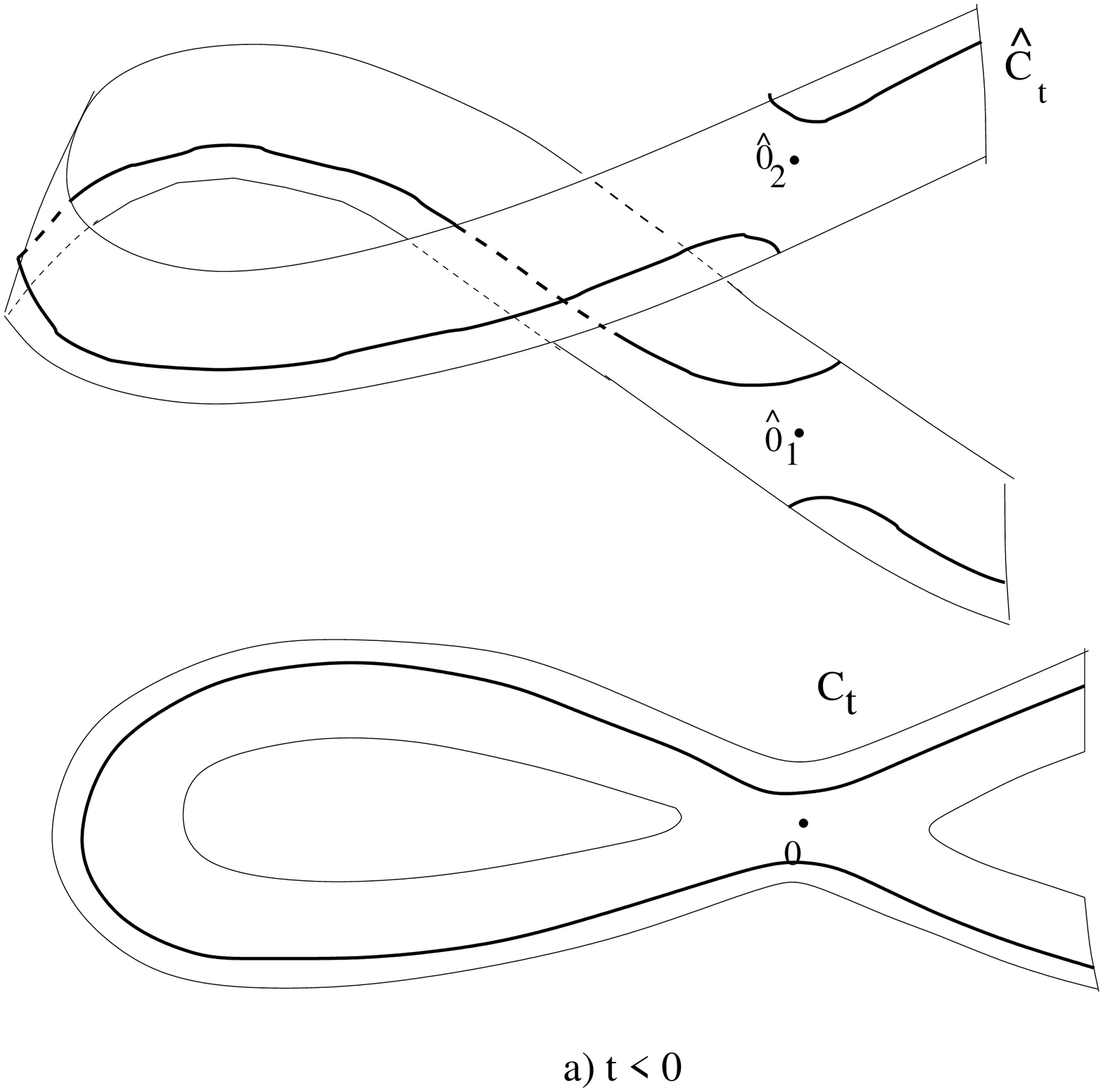}
\includegraphics[width=2.7in]{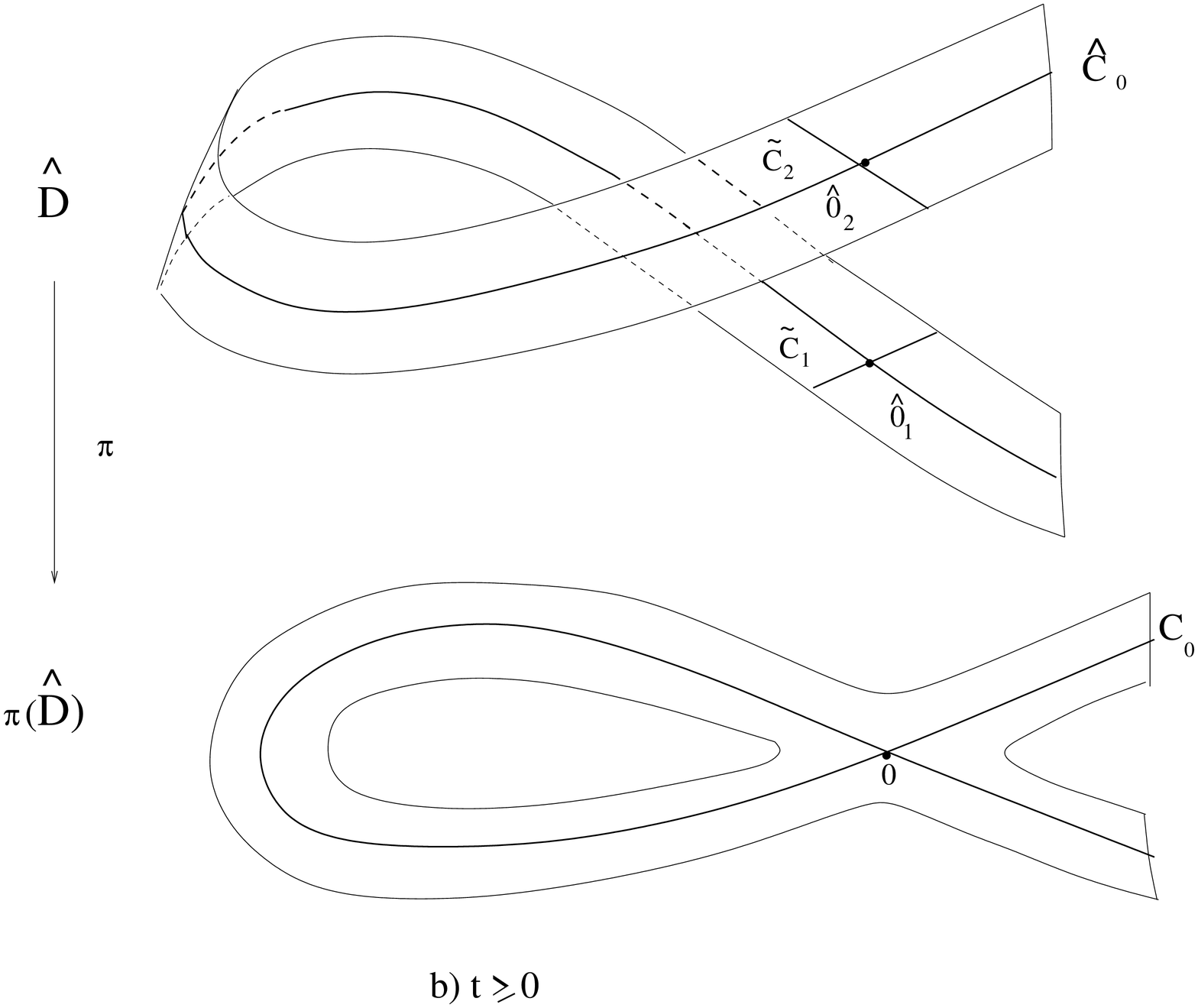}
\caption{This picture illustrates  Example \ref{moy-exmp1} and the proof of Theorem \ref{cpd}. 
On the right, for $t\ge 0$, $C_0\subset D$ has $0$ as a point of selfintersection. $\hat 0_1$ 
and $\hat 0_2$ are $\pi$-preimages of $0$ in $\widehat{D}$. $\widehat{C}_0$ is the lift of $C_0$ 
to the envelope and components $\tilde C_1$ and $\tilde C_2$ are components of $\widetilde{C}_0
\setminus \widehat{C}_0$. Curves $C_t$ for $t>0$ (not drawn there) behave similarly, only
without self-intersection. On the left we draw $C_t$ and $\widehat{C}_t$ for $t<0$. 
$\widehat{C}_t$ are irreducible but have more boundary components than just $\i(\d C_t)$.}
\label{band2}
\end{figure}

\begin{exmp} \rm 
\label{moy-exmp1}
Consider the following complex curve  in $\cc^2$
\[
C\deff \{(z,w)\in \cc^2: w^2 = z^3 + z^2\}.
\]
$C$ is immersed and has one point of self-intersection, the origin. The self-intersection of $C$ at 
zero is transverse. Indeed, $C$ can be parameterized as follows
\begin{equation}
\eqqno(param)
z = \lambda^2-1 \qquad \text{ and } \qquad w = \lambda(\lambda^2-1), \quad \lambda\in \cc .
\end{equation}
Denote by $\Phi_0(\lambda) = (\lambda^2-1, \lambda (\lambda^2-1))$ the parameterization map.  Then 
$\Phi_0(\pm 1)=0$ and 
\begin{equation}
\eqqno(trans)
\d_{\lambda}\Phi_0 (1) = (2\lambda , 3\lambda^2 -1)|_{\lambda = 1} = (2,2) \trans (-2,2) = 
\d_{\lambda}\Phi_0 (-1).
\end{equation}
Fix some $R\ge 10$ and take the disk $\Delta_R$ of radius $R$ in $\cc$. Set $C_0 \deff \Phi_0(\Delta_R)$. 
This is a compact analytic set with boundary $\d C_0 = \Phi_0(\d \Delta_R)$ in $\cc^2$. From \eqqref(param) 
we see that $w/z = \lambda$ is a holomorphic function on the normalization $C_0^{\n}$ of $C_0$. This function 
takes different values at two distinct points over the origin: $\lambda = \pm 1$. Imbed $\cc^2$ to $\cc^3$
as $\cc^2\times \{0\}$ and extend $\lambda $ as a multivalued analytic function $\Lambda$ to a 
$\delta$-neighborhood of $\bar C_0$ in $\cc^2$ for some $\delta>0$ and then consider $\Lambda$ as a function 
of three variables which doesn't depend on $u$. Consider the following family of parameterized 
compact curves with boundary in $\cc^3$
\begin{equation}
\eqqno(t>0)
C_t \deff 
\begin{cases}
z = \lambda^2 -1 \cr
w = \lambda(\lambda^2-1) \cr
u = t \lambda 
\end{cases}
\quad\text{where} \quad \lambda \in \bar\Delta_R \text{ and } 0 \le t \le \eps .
\end{equation}
Here $\eps >0$ is taken small enough in order that $|t\lambda|<\delta$ for $|t|\le \eps$ and $\lambda\in 
\bar\Delta$, \ie $\eps <\frac{\delta}{R}$. Notice that all $C_t$ are imbedded except for $C_0$.  Let 
$\Phi_t : \bar\Delta_R\to \cc^3$ 
\begin{equation}
\eqqno(phi)
\Phi_t: (\lambda ) \to (\lambda^2-1, \lambda (\lambda^2-1), t\lambda )
\end{equation}
for $t\in [0,\eps]$ be the parameterization of $C_t$. Consider the following family of convex cones 
``over $C_0$'', see the Picture \ref{cones}:
\begin{equation}
\eqqno(cones)
K^{\lambda}_{\delta} \deff \{\Phi_0(\lambda)\}\times \{u:\left|\arg u - \arg \lambda\right| < 
\arcsin \delta ,|u|<\delta \}.
\end{equation}
Point $\Phi_0(\lambda)$ let us call the {\slsf base point} of $K^{\lambda}_{\delta}$.
Remark that since $\Phi_0(1) = \Phi_0(-1)=0$ both 
\begin{equation}
\eqqno(cone1)
K^1_{\delta} = \{0\}\times \{u:\left|\arg u \right| < \arcsin \delta ,|u|<\delta \}
\end{equation}
and 
\begin{equation}
\eqqno(cone-1)
K^{-1}_{\delta} = \{0\} \times \{u: \left|\arg u -\pi \right| < \arcsin \delta ,|u|<\delta \}
\end{equation}
have the same base point - the origin. At the same time one observes that $K^1_{\delta}\cap 
K^{-1}_{\delta} = \emptyset$ provided $\delta >0$ was taken small enough.


\begin{claim}
\label{claim-4.1}
If $\eps > 0$ is taken small enough then for all $0\le t\le \eps$ one has 
\begin{equation}
\Phi_t (\lambda ) \in K^{\lambda}_{\delta} \quad\text{ for all } \quad 
\lambda\in \bar\Delta_R .
\end{equation}
\end{claim}

\noindent Indeed, the $(z,w)$-component of $\Phi_t (\lambda )$ belongs to $\bar C_0$, in fact it is 
$\Phi_0(\lambda)$. As for
the $u$-component we see that for $0<t<\eps $ we have that $\arg (t\lambda) = \arg (\lambda)$ and 
$|t\lambda |<\delta$ the latter satisfies the bound as in \eqqref(cones), see the Picture 
\ref{cones}. Set 
\begin{equation}
\eqqno(cones-f1)
\calk_{\delta} \deff \bigcup_{\lambda\in \bar\Delta _R}K_{\delta}^{\lambda}.
\end{equation}
We just proved that $C_k\subset \calk_{\delta}$ whatever small $\delta >0$ is provided that $k$ is 
big enough.

\begin{figure}[h]
\centering
\includegraphics[width=3.0in]{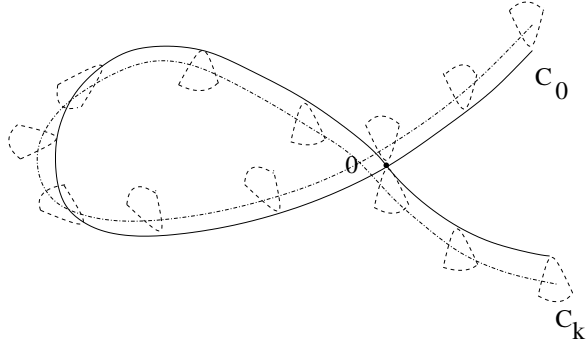}
\caption{Every $C_k$ for $k\gg 1$ enters to $\calk_{\delta}$. The latter is the union of the 
cones $K_{\delta}^{\lambda}$ with base points in $\bar C_0$ as on this picture.}
\label{cones}
\end{figure}

In order to make from $\calk_{\delta}$ a domain we need just to ``thicken'' it.  Using the
fact that the normal bundle to $C_0^{\n}$ is holomorphically trivial we can extend $\Phi_0$ to a
holomorphic immersion  
\begin{equation}
\widetilde{\Phi }_0: \bar\Delta_R\times \Delta_{\delta}  \to \cc^2.
\end{equation}
In $\Delta_{\delta}$ the variable we shall denote as $\mu$. And then thicken $\calk_{\delta}$ 
to a domain 
\begin{equation} 
\eqqno(dom-1)
D_{\delta} \deff \bigcup_{\lambda\in \bar\Delta _R, \mu\in \Delta_{\delta}}K_{\delta}^{\lambda ,\mu},
\end{equation}
where 
\begin{equation}
K_{\delta}^{\lambda ,\mu} \deff \{\Phi_0(\lambda , \mu )\}\times \{u:\left|\arg u - \arg \lambda\right| 
< \arcsin \delta ,|u|<\delta \}.
\end{equation}

We define $D$ as a union of two open sets $D_{\delta}$ and $D_2$, where  $D_2$ is a $\delta$-neighborhood 
of $\widetilde{\Phi}_0(\d\Delta_R\times \bar\Delta_{\delta})$. The latter is added to include 
the boundaries of $C_t$-s to $D$.

\smallskip Data $C_{\frac{1}{k}} \to C_0$ and 
$D$ satisfy the assumptions of Theorem \ref{cpd}. Restrict function $\Lambda$ to our $D$. This 
restriction is obviously singlevalued, but has only two-valued extension to a neighborhood of 
$0\in C_0$, \ie $\widehat{D}$ is two-sheeted over the origin and $C_0$ lifts to $\widehat{C}_0$ as 
on the Picture \ref{band2} b).

\medskip\noindent{\slsf Part II.} Now we shall extend our family $C_t$ to $t\in [-\eps ,0]$ providing the 
Example \ref{moy-exmp} from the Introduction. Set 
\begin{equation}
\eqqno(t<0)
C_t = \{(z,w,u): w^2 = z^3 + z^2 +t\}\cap \{\text{neighborhood of } C_0 \text{ in } \cc^2\}. 
\end{equation}
As such a neighborhood one can take $\widetilde{\Phi}_0(\Delta_R\times \Delta_{\delta})$ for example.
Notice that $C_t$ is contained in $\cc^2\times \{0\}$ for $t<0$ and is a continuous deformation of $C_0$
in the Hausdorff topology as well as in the (stronger) topology of currents. 
Our family $\{C_t\}_{t\in [-\eps,\eps]}$ is now constructed, \ie for $-\eps \le t\le 0$ curves $C_t$ are
defined as in \eqqref(t<0) and for $0\le t\le \eps$ as in \eqqref(t>0), see the Picture \ref{band2} below.
Notice that $\widehat{C}_t$ are close to $\widehat{C}_0\cup \widetilde{C}_1\cup\widetilde{C}_2$ for 
$t < 0$ and not to $\widehat{C}_0$ alone. This is why $\{\widehat{C}_t\}$ is discontinuous at 
$t=0$.
\end{exmp}

\begin{rema} \rm 
\label{ch-st}
We are bound at this point to examine the approach of \cite{CS}. It is claimed there roughly the 
following: {\it if a family holomorphic chains $\{C_t\}_{t\in [0,1]}$ is continuous in the topology 
of currents, $C_0\subset D$ and $\d C_t\subset D$ for all $t\in [0,1]$ then $C_1\subset \pi (\hat D)$. } 
As for the proof let $T$ be the set of such $t'$ that $C_t$ can be lifted to $\hat D$ up to $t'$. It is 
proved in \cite{CS} that $T$ is closed. This is equivalent to the ``discrete'' version of the CP in the 
topology of currents. After that it is claimed in \cite{CS} as an obvious fact that $T$ is also open. 
But our Example \ref{moy-exmp} shows that this is wrong. Therefore it seems natural to ask the following

\begin{problem} \rm 
Is the statement of \cite{CS} nevertheless holds true?
\end{problem}
May be one can ``modify'' the family and nevertheless achieve all points in $\hat D$ over $C_1$ by means 
of lifting this ``modified'' family?
\end{rema}

\newprg[LIFT-C.gr]{Gromov topology}

First let us recall few basic definitions concerning the Gromov topology on the space of complex curves
with boundary. For more details we refer to \cite{IS1,IS2}. Recall that a standard node is 
a complex analytic set $\calc_0 = \{(z_1,z_2)\in \Delta^2: z_1\cdot z_2=1\}$. A {\slsf nodal curve} is
a complex analytic set of pure dimension one with only nodes as singularities. A compact nodal curve with 
boundary $(C,\d C)$ is a nodal curve $C$ such that $\bar C$ is compact and smooth near its boundary 
$\d C \deff \bar C\setminus C$. Let $a_1,...,a_N$ be the nodes of $C$. We say that a real connected 
compact Riemann surface with boundary $(\Sigma , \d \Sigma)$ parameterizes $(C, \d C)$ if a continuous 
map $\sigma : \bar\Sigma \to \bar C$ is given such that: 

\smallskip

\sli  for every $k=1,...,N$ the set $\gamma_k\deff\sigma^{-1}(a_k)$ is a smooth imbedded circle in 
$\Sigma$;

\slii if $a_j\not= a_k$ then $\gamma_j\cap \gamma_k=\emptyset$;

\sliii $\sigma : \bar\Sigma \setminus \bigcup_{j=1}^N\gamma_k\to \bar C\setminus\{a_1,...,a_N\}$ is a 
diffeomorphism.

\smallskip By a complex curve  {\slsf over} a complex manifold $X$ we understand 
a pair $(C,u)$, where $C$ is a compact nodal curve with boundary and $u:C\to X$ is a holomorphic map 
continuous up to the boundary. Notice that the image $u(C)$ can have singularities other than just a
nodes.

\begin{defi} 
\label{gr-conv}
A sequence $(C_n,u_n)$ of complex curves over $X$ converges to a complex curve $(C_\infty,u_\infty)$ 
over $X$ in Gromov topology if all $C_n$ and $C_{\infty}$ can be parameterized 
by the same Riemann surface $\Sigma$ and the parameterizations $\sigma_n:\barr\Sigma \to \barr C_n$, 
$n\in \overline{\nn}\deff \nn\cup\{\infty\}$, can be chosen in such a way that the following hold:

\smallskip 

\sli $u_n\scirc \sigma_n$ converges to $u_\infty\scirc \sigma_\infty$ in the
$C^0(\bar\Sigma, X)$-topology, \ie uniformly on $\bar\Sigma$;

\slii if $\{ a_k \}$ is the set of nodes of $C_\infty$ and $\{\gamma_k\}$ are
the corresponding circles in $\Sigma$ then 

\quad on any compact subset $K\comp \bar\Sigma \bs \bigcup_k \gamma_k$ the convergence 
$u_n\scirc \sigma_n \to u_\infty \scirc \sigma_\infty$ is $\calc^{\infty}$;

\sliii for any compact subset $K\comp \barr\Sigma \bs \bigcup_k\gamma_k$ there
exists $n_0=n_0(K)$ such that 

\quad $ \sigma_n(K) \subset C_n\bs \{ nodes\} $
for all $n\geq n_0$ and the complex structures $\sigma_n^*j_{C_n}$ converge

\quad smoothly to $\sigma_{\infty}^*j_{C_{\infty}}$ on $K$;

\sliv the structures $\sigma_n^*j_{C_n}$ are constant in $n$ near the
boundary $\d\Sigma$.
\end{defi}

\smallskip Item (\sliv means that for every boundary circle $\gamma$ on $\Sigma$ there exist an 
annuli $A_n$ adjacent to $\sigma_n(\gamma)$ on $C_n$ (for all $n\in \overline{\nn}$) of a fixed 
conformal radius (\ie not depending on $n\in \overline{\nn}\}$). And tensors $\sigma_n^*j_{C_n}$ 
(where $j_{C_n}$ are tensors of complex structures on $A_n\subset C_n$) do not depend on $n\in 
\overline{\nn}\}$. This discussion leads us to the following

\begin{defi}
\label{gr-cont}
A family $\{(C_t,u_t)\}_{t\in [0,1]}$ of compact complex curves with boundary over $X$ is continuous 
in Gromov
topology if all $C_t$ can be parameterized by the same Riemann surface $\Sigma$ and parameterizations
$\sigma_t:\overline{\Sigma}_t \to \overline{C}_t$ can be chosen in such a way that the following hold:

\sli $u_t\circ \sigma_t$ is continuous as a function of a couple, \ie belongs to $\calc^0(\overline{\Sigma}
\times [0,1],X)$;

\slii for any $t_0\in [0,1]$ curve $(C_t,u_t)$ converge to $(C_{t_0}, u_{t_0})$ in the 
sense of Definition \ref{gr-conv} 

\quad when $t\to t_0$.
\end{defi}

Now let us see that the family $\{C_t\}$ of Example \ref{moy-exmp1} is discontinuous at zero in Gromov topology.
More accurately one should write $(C_t, \id )$, where $\id : C_t\to \cc^3$ is the natural inclusion, but we 
shall not do that. For $0\le t\le \eps$ curves $C_t$ are parameterized by the same $\Omega$ and parameterization
is given explicitly by $\Phi$ as in \eqqref(phi). So our family is continuous at zero from the right. 
For $-\eps \le t<0$ our family can be parameterized by torus with a hole. Indeed, for $-\eps \le t
\le 0$ all $C_t$, considered as algebraic curves in $\pp^2$, intersect the line at infinity at the same point 
$[0:1:0]$ and are smoothly imbedded, except of $C_0$ (provided $\eps >0$ was taken small enough). Remove the 
appropriate neighborhood of $[0:1:0]$ in $\pp^2$ and get toruses with a hole for all $-\eps \le t < 0$. When 
$t\nearrow 0$ an appropriate circle $\gamma_t$ on $C_t$ contracts to a point, thus producing a disk with one 
nodal point, which is $C_0$. So our family is continuous also from the left. But it cannot be parameterized by 
the same Riemann surface from the left and from the right because disc is not diffeomorphic to the torus with
a hole.

\newprg[LIFT-C.proof]{Lift to the envelope} 

Now we are going to prove Theorem \ref{cpc} from the Introduction. Denote by $T$ the set of $t'\in [0,1]$ such that 
$(C_t,u_t)$ can be continuously lifted to $\hat D$ up to $t'$. This means that for every $0\le t\le t'$ there exists 
a compact complex curve with boundary $(C_t,\hat u_t)$  over $\hat D$  such that 

\smallskip\sli $\hat u_t\circ\sigma_t:[0,t']\times \bar \Sigma \to \hat D$ is continuous as a mapping of two 
variables;

\slii $\pi \circ \hat u_t = u_t$ for all $t\in [0,t']$. 

\smallskip Let us underline that curves $C_t$ and their parameterizations $\sigma_t$ do not change. Notice also that 
(\slii implies that $\hat u_t = \i\circ u_t$ for $t$ close to zero. Our $T$ is obviously
non empty.

\smallskip\noindent{\slsf $T$ is closed.} Let $t_0\deff\supp\{t'\in T\}$.  As it was explained in Remark \ref{non-halph} 
$u_t(C_t)$ can be lifted to the envelope up to $t_0$ continuously in Hausdorff sense. Denote by $\widehat{C}_{t_0}$ 
the lift obtained this way. Take any $\hat y_0\in \widehat{C}_{t_0}$ and for $y_0\deff\pi (\hat y)$ consider a polydisk 
$\Delta^n(y_0,d)$ such that 
\begin{equation}
\eqqno(bihol)
\pi|_{\hat\Delta^n(\hat y_0,d)}:\hat\Delta^n(\hat y_0,d)\to \Delta^n(y_0,d)
\end{equation}
is a biholomorphism for an appropriate neighborhood $\hat\Delta^n(\hat y_0,d)$ of $\hat y_0$, see Remark 
\ref{bihol-pi}. Let $\hat x_0\in C_{t_0}$ be such that $u_{t_0}(\hat x_0)= y_0$ and $x_0\in \bar\Sigma$ such that 
$\sigma_{t_0}(x_0) = \hat x_0$. It may happen that several branches of $u_{t_0}(C_{t_0})$ pass 
through $y_0$, we fix one of them. Now we can set 
\begin{equation}
\eqqno(lift-u1)
\hat u_{t_0}(\hat x) \deff (\pi^{-1}\circ u_{t_0})(\hat x)
\end{equation}
for $\hat x$ close $\hat x_0$ on $C_{t_0}$ on this branch. This is well defined and holomorphic near $\hat x_0$.
We do this for all branches of $u_{t_0}(C_{t_0})$ passing through $y_0$. Performing this in a neighborhood of 
every $\hat y_0\in \widehat{C}_{t_0}$ we obtain a holomorphic lift $\hat u_{t_0}$ of $u_{t_0}$ as required.

\begin{rema} \rm In order to obtain the limit $\widehat{C}_{t_0}$ one may apply also the Gromov compactness
theorem in the form that is proved in \cite{IS2}. Note that $\hat D$ is holomorphically convex and therefore
all $\widehat{C}_t$ do stay in a compact part of $\hat D$ as $t\nearrow t_0$. But we don't need to use such 
strong statement here.
\end{rema}

\begin{figure}[h]
\centering
\includegraphics[width=5in]{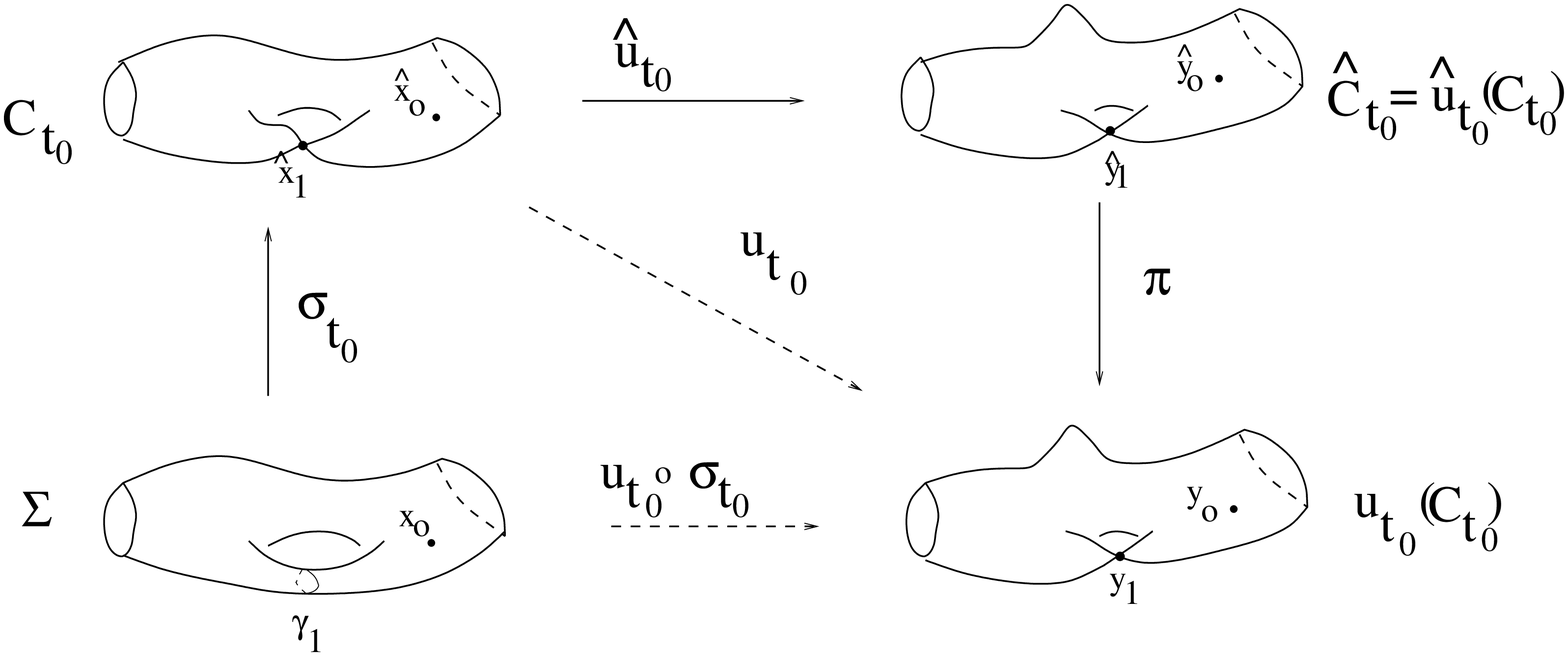}
\caption{This picture explains the definition of $\phi$ and then of $\hat u$ for $\hat x$ being a node (like $\hat x_1$)
or not (like $\hat x_0$).}
\label{lift}
\end{figure}

\smallskip\noindent{\slsf $T$ is open.}  We obviously have that $\widehat{C}_{t_0}\deff \hat u_{t_0}(C_{t_0})$. 
Take again any $\hat y_0\in \widehat{C}_{t_0}$ and for $y_0\deff\pi (\hat y_0)$ consider polydisks $\Delta^n(y_0,d)$ 
and $\hat\Delta^n(\hat y_0,d)$ as in \eqqref(bihol). Let $x_0\in \Sigma$ be such that 
$(u_{t_0}\circ\sigma_{t_0})(x_0) =y_0$ and let $\hat x_0\deff \sigma_{t_0}(x_0)\in C_{t_0}$. Define a mapping $\phi$ in a 
neighborhood $U$ of $(t_0,x_0)$ in $[0,1]\times \bar\Sigma$ with values in $\hat D$ as follows
\begin{equation}
\eqqno(lift-u2)
\phi :(t,x) \to (\pi^{-1}\circ u_t\circ \sigma_t)(x).
\end{equation}
$\phi$ is well defined and satisfies $(\pi\circ \phi)(t,x) = (u_t\circ \sigma_t)(x)$ for all $(t,x)$ in its domain of definition.

If $\hat x_0$ is not a node we can set $\hat u_t = \phi\circ \sigma_t^{-1}$. This map is well defined holomorphic on an appropriate 
open set (more exactly on $\sigma_t(U_t)$, where $U_t\deff U\cap \{t\}\times \bar\Sigma$) of every $C_t$ minus nodes (if any). At the 
same time it is clearly extends continuously to nodal points. Therefore it is holomorhic everywhere. Theorem is proved.

\smallskip\qed

\begin{exmp} \rm
\label{r2-remov}
To illustrate the situation in this theorem consider the following (known) example. Let 
$D=\cc^2\setminus \rr^2$ be a domain in $X=\cc^2$. Consider the following family of complex 
curves with boundary 
\[
C_t =\{(z_1+iz_2)(z_1-iz_2) = t\}\cap \bb^2(1), \quad\text{ where }\quad t\in [-1/2,1/2].
\]
One easily checks that $C_t\subset D$ for $t\in ]0,1/2]$ and boundaries $\d C_t$ stay in the compact 
\[
K=\{ z=x+iy: 1/4 \le \norm{x}^2\le 3/4, 1/4\le \norm{y}^2\le 3/4\} \comp D \text{ for {\slsf all} } t\in [-1/2,1/2].
\]
Finally as $t\to 0$ curves $C_t$ degenerate to a node $C_0 = \{z_1+iz_2=0\}\cup \{z_1-iz_2=0\}$,
and the latter contains the origin. This proves, via Theorem \ref{cpc}, the (well known) fact that 
$\rr^2$ is removable singularity for holomorphic/meromorphic functions of two variables.
\end{exmp}

\smallskip Finally we give the construction of the  Example \ref{nash-exmp} from the Introduction
showing that the assumption of Steiness of the ambient manifold $X$ cannot be dropped neither in 
Theorem \ref{cpd} nor in Theorem \ref{cpc}. This example was communicated to me by Chirka and first 
published in \cite{IS1}.

\begin{exmp} \rm
\label{r2-exmp}
As a complex manifold $X$ take the total space of the holomorphic rank two bundle $\calo(-1)\oplus \calo (-1)$
over the Riemann sphere $\pp^1$. By $z$ denote the standard affine coordinate on $\pp^1$. Let 
$\eta_1 = z\xi_1$ and $\eta_2 = z\xi_2$ be the standard coordinates on fibers. Denote by  $E = \{\xi_1=\xi_2=0\}$ 
the zero section of this bundle and by $\Sigma_1$ the complex hypersurface $\{\xi_1=0\}$. Consider the 
following {\slsf holomorphic} function $f = e^{\xi_2 /\xi_1}$ on $D=X\setminus \Sigma_1$. Remark that $\Sigma_1$ in
an essential singularity of $f$, \ie $f$ doesn't extend to a neighborhood of any point of $\Sigma$ even
meromorphically.

\smallskip Let $[z_1:z_2]$ be the homogeneous coordinates in $\pp^1$, and $z=z_2/z_1$ is our affine coordinate. Let us
restrict ourselves to the hypersurface $\Sigma_2 = \{\xi_2=0\}$ where our curves $C_t$ will live. This hypersurface is 
the blown up $\pp^2$ and $z_1,z_2$ can be considered as the affine coordinates with $\xi_1=z_1$ and $\eta_1=z_2$. Indeed 
$z\xi_1 = z_2/z_1\cdot z_1 = z_2= \eta_1$ as it should be. In these affine coordinates we consider the following family
of complex curves (analytic disks in fact) with boundary: $C_t \deff \{ |z_1|\le 1, z_2 = t, \xi_2 = 0 \} \subset 
\Sigma_2\cap D$, $t\in [0,1]$. The limit of this family as $t\searrow 0$ is $C_0 = E\cup \bar\Delta $, where $\bar\Delta
\deff \{ |z_1|\le 1, z_2=0, \eta_2=0\}$. We are well under the conditions of Theorem \ref{cpc} except of non-Steiness of $X$. 
But we see that our holomorphic in $D$ function $f$ has an essential singularity on $E\subset C_0$.

\smallskip If we take $C_{\frac{1}{k}}$ as $C_k$ we get a counterexample to Theorem \ref{cpd} for non-Stein $X$.
\end{exmp}

\smallskip  One can remark that our Examples \ref{moy-exmp} and \ref{nash-exmp} do live in dimension three. On our opinion in 
dimension two the continuous version of CP (as well as a discrete one) holds true. Note that the Cartan-Thullen construction 
provides us with an envelope of meromorphy $(\hat D, \pi)$ any domain $D$ in any complex manifold.

\begin{problem}
Let $D$ be a domain in a complex surface $X$ and let $\{(C_t,u_t)\}_{t\in [0,1]}$ be a continuous in Gromov topology 
family of stable curves over $X$ such that:

\sli $u_0(C_0)\subset D$;

\slii $u_t(\d C_t)\subset D$ for all $t\in [0,1]$.

\noindent Then $\{(C_t,u_t)\}$ can be lifted to the envelope of meromorphy $\hat D$ of $D$ as in Theorem \ref{cpc}. 
\end{problem}

A weaker statement was proved in \cite{IS1}, see Theorems 2.2.2 and 2.2.3 there.

\ifx\undefined\bysame
\newcommand{\bysame}{\leavevmode\hbox to3em{\hrulefill}\,}
\fi

\def\entry#1#2#3#4\par{\bibitem[#1]{#1}
{\textsc{#2 }}{\sl{#3} }#4\par\vskip2pt}

\end{document}